\documentclass[leqno,11pt]{article}

\usepackage{amssymb,amsmath,amscd,amsfonts,a4,psfrag,graphicx,latexsym,epsfig,color}
\usepackage[all]{xy}

\usepackage{framed}

\usepackage{tikz-cd}
\usetikzlibrary{cd}

\usepackage{verbatim}

\newtheorem{theorem}{Theorem}[section]
\newtheorem{lemma}[theorem]{Lemma}
\newtheorem{proposition}[theorem]{Proposition}
\newtheorem{corollary}[theorem]{Corollary}
\newtheorem{definition}[theorem]{Definition}
\newtheorem{remark}[theorem]{Remark}

\newtheorem{assumption}[theorem]{Assumption}

\newcommand{\ggot}{\ensuremath{\mathfrak{g}}}

\newcommand{\kgot}{\ensuremath{\mathfrak{k}}}

\newcommand{\pgot}{\ensuremath{\mathfrak{p}}}

\newcommand{\zgot}{\ensuremath{\mathfrak{z}}}

\newcommand{\tgot}{\ensuremath{\mathfrak{t}}}
\newcommand{\ugot}{\ensuremath{\mathfrak{u}}}
\newcommand{\Rgot}{\ensuremath{\mathfrak{R}}}
\newcommand{\Sgot}{\ensuremath{\mathfrak{S}}}
\newcommand{\Xgot}{\ensuremath{\mathfrak{X}}}

\newcommand{\glgot}{\ensuremath{\mathfrak{gl}}}

\newcommand{\Acal}{\ensuremath{\mathcal{A}}}

\newcommand{\Fcal}{\ensuremath{\mathcal{F}}}

\newcommand{\Ocal}{\ensuremath{\mathcal{O}}}
\newcommand{\Pcal}{\ensuremath{\mathcal{P}}}

\newcommand{\Scal}{\ensuremath{\mathcal{S}}}
\newcommand{\Tcal}{\ensuremath{\mathcal{T}}}

\newcommand{\Z}{\ensuremath{\mathbb{Z}}}
\newcommand{\C}{\ensuremath{\mathbb{C}}}

\newcommand{\G}{\ensuremath{\mathbb{G}}}

\newcommand{\Q}{\ensuremath{\mathbb{Q}}}
\newcommand{\R}{\ensuremath{\mathbb{R}}}

\newcommand{\croc}{\ensuremath{\hookrightarrow}}

\newcommand{\diag}{{\ensuremath{\rm Diag}}}

\def \tG {\widetilde{G}}
\def \tK {\widetilde{K}}
\def \tU {\widetilde{U}}
\def \tB {\widetilde{B}}
\def \tP {\widetilde{P}}

\def \tW {\widetilde{W}}

\def \tOcal {\widetilde{\Ocal}}

\def \tggot {\tilde{\ggot}}

\def \tkgot {\tilde{\kgot}}
\def \tpgot {\tilde{\pgot}}
\def \tugot {\tilde{\ugot}}
\def \ttgot {\tilde{\tgot}}

\def \tX {\widetilde{X}}

\def \tT {\widetilde{T}}
\def \tP {\widetilde{P}}

\def \txi {\tilde{\xi}}

\def \tw {\tilde{w}}

\newcommand{\e}{\ensuremath{\mathrm{e}}}

\newcommand{\s}{\ensuremath{\mathrm{s}}}

\newcommand{\horn}{\ensuremath{\hbox{\rm Horn}}}

\newcommand{\LR}{\ensuremath{\hbox{\rm LR}}}

\begin{document}

\title{Eigenvalues, singular values, and the O'Shea-Sjamaar Theorem}

\author{Paul-Emile Paradan\footnote{IMAG, Univ Montpellier, CNRS, email : paul-emile.paradan@umontpellier.fr}}

\maketitle

\date{}

%%%%%%%%%%%%%%%%%%%%%%%%%%%%%%%%%%%%%%%%%%%%%%%%%%%%%%%%%%%%%%%%%%%%%

\begin{abstract}
The main focus of this work is the study of several cones relating the eigenvalues or singular values of a matrix to those of its off-diagonal blocks.
\end{abstract}

%%%%%%%%%%%%%%%%%%%%%%%%%%%%%%%%%%%%%%%%%%%%%%%%%%%%%%%%%%%%%%%%%%%%%

%{\def\thefootnote{\relax}
%\footnote{{\em Keywords} : 
%transversally elliptic symbol.\\
%{\em 1991 Mathematics Subject Classification} : 58F06, 57S15, 19L47.}
%\addtocounter{footnote}{-1}
%}

\tableofcontents

%%%%%%%%%%%%%%%
%%%%%%%%%%%%%%%
\section{Introduction}
%%%%%%%%%%%%%%%
%%%%%%%%%%%%%%%

Let $p\geq q\geq 1$ and $n=p+q$. Let $Herm(n)$ denote the vector space of $n$-square Hermitian matrices. The spectrum of $X\in Herm(n)$ is denoted by 
$\e(X)=(\e_1\geq \cdots\geq\e_n)$ and the singular spectrum of 
a matrix $Y\in M_{p,q}(\C)$ is denoted by $\s(Y)=(\s_1\geq \cdots \geq \s_q\geq 0)$.

The main purpose of this article is to describe the following cones:
\begin{align*} 
\Acal(p,q) &= \Big\{(\e(X), \s(X_{12})), \ X\in Herm(n)\Big\}, \\ 
\Scal(p,q) &= \Big\{(\s(X), \s(X_{12}),\s(X_{21})), \ X\in M_{n,n}(\C)\Big\}, \\
\Tcal(p,q) &= \Big\{(\s(X), \s(X_{11}),\s(X_{22})), \ X\in M_{n,n}(\C)\Big\}.
\end{align*}
Here, a $n$-square complex matrix $X$ is written by blocks 
$X=\begin{pmatrix}
X_{11}& X_{12}\\
X_{21}& X_{22}
\end{pmatrix}$
where $X_{12}\in M_{p,q}(\C)$ and $X_{21}\in M_{q,p}(\C)$.

In the 1970s, Thompson gave some inequalities satisfied by the elements of $\Tcal(p,q)$ \cite{Thompson72,Thompson75}, and more recently 
Fomin, Fulton, Li and Poon obtained sets of inequalities that describe the cone $\Acal(p,q)$ \cite{Li-Poon04,FFLP}.

\medskip

The main objective of this work is to explain how a direct application of O'Shea-Sjamaar's theorem \cite{OSS} yields complete sets of inequalities 
for the cones $\Acal(p,q)$, $\Scal(p,q)$, and $\Tcal(p,q)$. However, this method does not provide an optimal description of these cones, as 
it leads to a large number of redundancies in the list of inequalities.  We'll see, for example, that the Fomin-Fulton-Li-Poon description of 
$\Acal(p,q)$ is more accurate than ours. In a future work \cite{pep-real-3}, we will propose a more precise method to describe these inequality sets, using the main result of \cite{pep-real-2}.

\medskip

Throughout this article, we make extensive use of Horn cones $\horn(n)$ and Littlewood-Richardson cones $\LR(m,n)$. Let's recall their definition.
To any integers $n,m\geq 1$, we associate
\begin{eqnarray*}
\horn(n)&=& \Big\{(\e(X), \e(Y),\e(X+Y)), \ X,Y\in Herm(n)\Big\}\\
\LR(m,n)&=&\Big\{(\e(M), \e(M_{\bf I}),\e(M_{\bf II})), \ M\in Herm(m+n)\Big\},
\end{eqnarray*}
where $M_{\bf I}\in Herm(m)$ and $M_{\bf II}\in Herm(n)$ are the extracted matrices such that 
$M=\begin{pmatrix}
M_{\bf I}& *\\
*& M_{\bf II}
\end{pmatrix}$.

In \S \ref{sec:classical}, we recall the description obtained by Klyachko \cite{Klyachko} and Knutson-Tao \cite{Knutson-Tao-99} for $\horn(n)$ cones, and 
 that obtained by Berenstein-Sjamaar \cite{Berenstein-Sjamaar} and Ressayre \cite{Ressayre10} of $\LR(m,n)$ cones. In both cases, 
 the inequalities are parameterized using Littlewood-Richardson coefficients.

In \S \ref{sec:A-p-q}, we show that $\Acal(p,q)$ can be characterized as a sub-cone of\footnote{Here $n=p+q$.} $\horn(n)$. To any $\lambda=(\lambda_1\geq \cdots\geq \lambda_n)$ and 
$s=(s\geq \cdots\geq s_q\geq 0)$, we associate $\lambda^*=(-\lambda_n\geq \cdots\geq -\lambda_1)$ and 
$$
\widehat{s}\,^{p,q}:=(s_1\geq\cdots\geq s_q\geq\underbrace{0,\cdots,0}_{p-q}\geq-s_q\geq\cdots\geq-s_1).
$$
We then show that  $(\lambda,s)\in \Acal(p,q)$ if and only if $(\lambda,\lambda^*, 2\widehat{s}\,^{p,q})\in \horn(n)$. This makes it possible to describe $\Acal(p,q)$ by means of the inequalities defining $\horn(n)$, but we'll see that the resulting description is less precise than that given by Fomin-Fulton-Li-Poon in \cite{FFLP}.

In \cite{FFLP}, the authors pose the question of finding a collection of linear inequalities that describes $\Scal(p,q)$ (Problem 1.15). 
We answer this problem in \S \ref{sec:S(p,q)} by showing that $\Scal(p,q)$ can be characterized as the intersection of $\LR(n,n)$ with the subspace of 
$\R^{2n}\times\R^{n}\times\R^n$ formed by the elements $(\widehat{\gamma}^{n,n},\widehat{s}\,^{p,q}, \widehat{t}\,^{p,q})$ where $(\gamma,s,t)\in \R^{n}\times\R^{q}\times\R^q$.

In the last section, we give a set of inequalities describing the cone $\Tcal(p,q)$, showing that $\Tcal(p,q)$ is characterized as the intersection of 
$\LR(2p,2q)$ with the subspace of $\R^{2n}\times\R^{2p}\times\R^{2q}$ formed by the elements $(\widehat{\gamma}^{n,n},\widehat{s}\,^{p,p}, \widehat{t}\,^{q,q})$ 
where $(\gamma,s,t)\in \R^{n}\times\R^{p}\times\R^q$.

In the case $\Tcal(p,1)$, we recover the interleaving inequalities of singular values obtained by Thompson \cite{Thompson72}.

%%%%%%%%%%%%%%%%%%%%%
\section*{Acknowledgement}
%%%%%%%%%%%%%%%%%%%%%

I would like to thank Nicolas Ressayre for the interesting discussions we had on this subject. I'm also grateful to a group of young students, Martina Ag\"{u}era Sanchez and  Ariane $\&$ Constantin Paradan, who implemented A.S. Buch's ``Littlewood-Richardson calculator" in various programs, enabling me to calculate a few examples. I would also like to thank the referees for their comments, which helped me improve this text.

%%%%%%%%%%%%%%%%%%%%%
\section*{Notations}
%%%%%%%%%%%%%%%%%%%%%

Throughout the paper :
\begin{itemize}
\item We fix $p\geq q\geq 1$ and $n=p+q$. 
\item We write $0_{ab}$ for the zero matrix of size $a\times b$.
\item Let $M_{a,b}(\C)$ be the vector space of complex $a \times b$ matrices.
\item $\R^\ell_+$ is the set of sequences $x=(x_1\geq \cdots \geq x_\ell)$ of real numbers.
\item $\R^\ell_{++}$ is the set of sequences $x=(x_1\geq \cdots \geq x_\ell\geq 0)$ of non-negative real numbers.
\item For any positive integer $\ell$, let $[\ell]$ be the set $\{1,\ldots,\ell\}$.
\item If $x\in\R^\ell$ and $A\subset [\ell]$, we write $|x|_A=\sum_{a\in A}x_a$ and $|x|=\sum_{i=1}^\ell x_i$.
\item For $A\subset [\ell]$, we define $A^o:=\{\ell+1- a, a\in A\}$ and $A^c:=[\ell]\setminus A$.
\item If $x\in\R^\ell$, let $\diag(x)$ be the diagonal $\ell\times\ell$ matrix with diagonal entries equal to $x_1,\ldots,x_\ell$.
\item If $A=\{a_1<\cdots< a_p\}$ is an increasing sequence of positive integers, let $\mu(A)=(a_p-p\geq \cdots \geq a_1-1\geq 0)$.
\end{itemize}

%%%%%%%%%%%%%%%%%%%%%%%%%%%%%%%%%%%%%%%%%%%%%%%%
%%%%%%%%%%%%%%%%%%%%%%%%%%%%%%%%%%%%%%%%%%%%%%%%
\section{Reminder of some classical results}\label{sec:classical}
%%%%%%%%%%%%%%%%%%%%%%%%%%%%%%%%%%%%%%%%%%%%%%%%
%%%%%%%%%%%%%%%%%%%%%%%%%%%%%%%%%%%%%%%%%%%%%%%%

We recall some classical facts that we'll be needing later on.

%%%%%%%%%%%%%%%%%%%%%%%%%%%%%
\subsection{Singular values}\label{sec:singular-values}
%%%%%%%%%%%%%%%%%%%%%%%%%%%%%

Let $X$ be a rectangular matrix, say $m\times n$, with complex entries, and let $X^*$ denote the complex conjugate transpose of $X$. 
Let $\eta_1(X)\geq \cdots\geq \eta_m(X)\geq 0$ be the eigenvalues of the positive semidefinite matrix $XX^*$. Notice that $\eta_k(X)=0$ 
when $k>\ell:=\inf\{m,n\}$.

The singular values of the matrix $X$ are the coordinates of the vector 
$$
\s(X):=\left(\sqrt{\eta_1(X)},\ldots,\sqrt{\eta_\ell(X)}\right)\in \R^\ell_{++}.
$$

Consider the canonical action of the unitary group $U_m\times U_n$ on $M_{m,n}(\C)$: $(g,h)\cdot X= gXh^{-1}$, $\forall (g,h)\in U_m\times U_n$. 
The singular values map  $\s : M_{m,n}(\C)\to \R^\ell_{++}$ induces a bijective application  \ 
$M_{m,n}(\C)/U_m\times U_n\stackrel{\sim}{\longrightarrow}\R^\ell_{++}$.

%%%%%%%%%%%%%%%%%%%%%%%%%%%%%
\subsection{Augmented matrices}\label{sec:augmented}
%%%%%%%%%%%%%%%%%%%%%%%%%%%%%

Let $p\geq q\geq 1$. 

If $Y$ is a $p\times q$ matrix, we denote by $\widehat{Y}^{p,q}$ the $n$-square Hermitian matrix $\begin{pmatrix}0_{pp}& Y\\Y^*& 0_{qq}\end{pmatrix}$. Here 
the spectrum of $\widehat{Y}^{p,q}$ is equal to $\widehat{s}\,^{p,q}:=(s_1,\cdots,s_q,0,\cdots,0,-s_q,\cdots,-s_1)\in \R^n_+$, where $s\in \R^q_{++}$ is the singular spectrum of $Y$.

If $Z$ is a $q\times p$ matrix, we denote by $\widehat{Y}^{q,p}$ the $n$-square Hermitian matrix $\begin{pmatrix}0_{qq}& Z\\Z^*& 0_{pp}\end{pmatrix}$. 
Here the spectrum of $\widehat{Z}^{q,p}$ is also equal to $\widehat{t}\,^{p,q}$, where $t=\s(Z)\in \R^q_{++}$.

If $X$ is a $k$-square matrix, we simply denote $\widehat{X}^{k,k}$ by $\widehat{X}$. The spectrum of $\widehat{X}$ is also simply denoted 
$\widehat{\mu}:=(\mu_1,\cdots,\mu_k,-\mu_k,\cdots,-\mu_1)\in \R^{2k}_+$,  where $\mu=\s(X)$.

%%%%%%%%%%%%%%%%%%%%%%%%%%%%%
\subsection{Horn inequalities }
%%%%%%%%%%%%%%%%%%%%%%%%%%%%%
 
%Recall that $\e(X)\in\R^n_+$ denotes the vector of eigenvalues of a $n$-square Hermitian matrix $X$. We consider the Horn cone 
%$\horn(n):=\left\{(\e(X),\e(Y),\e(X+Y));\ X,Y\in Herm(n) \right\}$. 

Denote the set of cardinality $r$-subsets $I=\{i_1<i_2<\cdots<i_r\}$ of $[n]$ by $\Pcal^n_r$.

\begin{definition}
For any $1\leq r<n$, $\LR^n_r$ refers to the set of triplets $(I,J,K)\in (\Pcal^n_r)^3$ such that $(\mu(I),\mu(J),\mu(K))\in\horn(r)$.
\end{definition}

The following theorem was conjectured by Horn \cite{Horn} and proved by a combination of the works of Klyachko \cite{Klyachko} and Knutson-Tao 
\cite{Knutson-Tao-99}. 

\medskip

\begin{theorem}\label{theo:horn} The triplet $(x,y,z)\in (\R^n_+)^3$ belongs to $\horn(n)$ if and only if the following conditions hold:
\begin{itemize}
\item $|x|+|y|=|z|$,
\item $|x|_I+|y|_J\geq |z|_K$, for any $r<n$ and any $(I,J,K)\in \LR^n_r$.
\end{itemize} 
\end{theorem}

In the following sections, we'll use Littlewood-Richardson coefficients to parameterize certain inequalities. Let's recall their definition. 
Let $\lambda,\mu$, and $\nu$ be three partitions of length less than $n\geq 1$. We associate them with the irreducible representations $V_\lambda$, $V_\mu$ and $V_\nu$ 
of the unitary group $U_n$. The Littlewood-Richardson coefficient $c^\lambda_{\mu,\nu}$ can be characterized by the relation
$c^\lambda_{\mu,\nu}=\dim \left[V_\lambda^*\otimes V_\mu\otimes V_\nu\right]^{U_n}$. Thanks to the saturation Theorem of Knutson and Tao 
\cite{Knutson-Tao-99}, we know that
$c^\lambda_{\mu,\nu}\neq 0\Longleftrightarrow (\mu,\nu,\lambda)\in \horn(n)$.

The following kind of duality is used in the next sections: for all $(I,J,K)\in (\Pcal^n_r)^3$, we have
\begin{equation}\label{eq:duality}
c^{\mu(K)}_{\mu(I),\mu(J)}\neq 0\Longleftrightarrow c^{\mu((K^o)^c)}_{\mu((I^o)^c),\mu((J^o)^c)}\neq 0.
\end{equation}
Since the relation 
$|x|_{(I^o)^c}+ |y|_{(J^o)^c} \geq |z|_{(K^o)^c}$ is equivalent to $|x|_{I^o}+|y|_{J^o} \leq |z|_{K^o}$, in Theorem \ref{theo:horn}, we can rewrite  the last condition by requiring that 
$$
|x|_I+|y|_J\geq |z|_K\quad {\rm and}\quad|x|_{I^o}+|y|_{J^o}\leq |z|_{K^o}
$$
for any $r\leq \tfrac{n}{2}$ and any $(I,J,K)\in \LR^n_r$.

%%%%%%%%%%%%%%%%%%%%%%%%%%%%%
\subsection{The cone $\LR(U,\tU)$}
%%%%%%%%%%%%%%%%%%%%%%%%%%%%%

Let $\iota:U\croc \tU$ be two connected compact Lie groups. We choose an invariant scalar product $(-,-)$ on the Lie algebra $\tugot$ of $\tU$, and we denote by $\pi:\tugot\to\ugot$ 
the orthogonal projection.

Select maximal tori $T$ in $U$ and $\tT$ in $\tU$ such that $T\subset \tT$, and
Weyl chambers $\tgot_+\subset\tgot$ and $\ttgot_+\subset\ttgot$, where $\tgot$ and $\ttgot$ denote the Lie algebras of $T$, resp. $\tT$. 
The aim of this section is to recall the description of the following cone given in \cite{Berenstein-Sjamaar,Ressayre10}:
$$
\LR(U,\tU)=\left\{(\xi,\txi)\in \tgot_+\times\ttgot_+,\ U\xi\subset \pi\big(\,\tU\txi\,\big)\right\}.
$$

Consider the lattice $\wedge:=\frac{1}{2\pi}\ker(\exp : \tgot\to T)$ and the Weyl groups $\tW=N_{\tU}(\tT)/\tT$ and $W=N_{U}(T)/T$. We denote by $w_o\in W$ the longest element. 
A vector $\gamma\in \tgot$ is called {\em rational} if it belongs to the $\Q$-vector space $\tgot_\Q$ generated by $\wedge$. We will see that the cone 
$\LR(U,\tU)$ is completely described by inequalities of the form 
$$
(\tilde{\xi},\tilde{w}\gamma)\geq (\xi,w_ow\gamma)
$$
with $\gamma$ rational anti-dominant and $(w,\tw)\in W\times \tW$.

%%%%%%%%%%%%%%%%%%%%%%%%%%%
\subsubsection{Admissible elements}
%%%%%%%%%%%%%%%%%%%%%%%%%%%

We let $\Sigma(\tugot/\ugot)\subset \tgot^*$ denote the set of weights relative to the $T$-action on $(\tugot/\ugot)\otimes\C$. If $\gamma\in\tgot$, we denote by 
$\Sigma(\tugot/\ugot)\cap \gamma^\perp$ the subset of weights vanishing against $\gamma$.

\begin{definition}\label{admissible-exemple-1} A rational element $\gamma\in \tgot$ is said {\em admissible} when
\begin{equation}\label{eq:condition-gamma}
{\rm Vect}\big(\Sigma(\tugot/\ugot)\cap \gamma^\perp\big)={\rm Vect}\big(\Sigma(\tugot/\ugot)\big)\cap \gamma^\perp.
\end{equation}
\end{definition}

\medskip 

When the following assumption is satisfied, we'll see in Section \S \ref{sec:LR-U-tildeU} that $\LR(U,\tU)$ is described by  inequalities parameterized by a finite number of admissible elements.

\begin{assumption}\label{hypothese}
The subspace $\zgot:=\{X\in \tgot, \alpha(X)=0, \forall \alpha\in\Sigma(\tugot/\ugot)\}$ is contained in the center  $Z_{\tugot}$ of $\tugot$. 
\end{assumption}
This assumption means that any ideal of $\tugot$ contained in $\ugot$ is a subspace of $Z_{\tugot}\cap\ugot$.

Let $\tgot=\zgot\oplus\tgot_1$ be a rational decomposition. Let us denote by $\Sigma(\tugot/\ugot)'$ the image of $\Sigma(\tugot/\ugot)$ 
through the projection $\tgot^*\to (\tgot_1)^*$. If Assumption \ref{hypothese} holds,  $\Sigma(\tugot/\ugot)'$ generates $(\tgot_1)^*$. Any rational element $\gamma\in\tgot$ can be written 
$\gamma=\gamma_0+\gamma_1$ where $\gamma_0$ is a rational element of $\zgot$ and $\gamma_1$ is a rational element of $\tgot_1$. 
We see then that a rational element $\gamma$ is admissible if and only if $\gamma_1$ is admissible. The later condition is equivalent to asking that 
the hyperplane $(\gamma_1)^\perp\subset (\tgot_1)^*$ is generated by a finite subset of $\Sigma(\tugot/\ugot)'$. Thus, there are a finite number of choices for $\gamma_1$ 
(up to multiplication by $\Q^{>0}$).

%%%%%%%%%%%%%%%%%%%%%%%%%%%
\subsubsection{Polarized trace}
%%%%%%%%%%%%%%%%%%%%%%%%%%%

Let $\Rgot(\ugot)$ and $\Rgot(\tugot)$ be the set of roots associated to the Lie algebras $\ugot$ and $\tugot$. The choice of the Weyl chambers $\tgot_{+}$
and $\ttgot_{+}$  define subsets of positive roots $\Rgot^+\subset\Rgot(\ugot)$ and $\widetilde{\Rgot}^+\subset\Rgot(\tugot)$.

For a rational element $\gamma\in\tgot$ and $(w,\tw)\in W\times \tW$, we will use the following condition to parameterize the inequalities of $\LR(U,\tU)$:

\begin{equation}\label{eq:trace-condition-w-tilde}
\sum_{\stackrel{\alpha\in\Rgot^+}{\langle\alpha,w\gamma\rangle> 0}}\langle\alpha,w\gamma\rangle+
\sum_{\stackrel{\tilde{\alpha}\in\widetilde{\Rgot}^+}{\langle\tilde{\alpha},\tilde{w}\gamma\rangle< 0}}\langle\tilde{\alpha},\tilde{w}\gamma\rangle=0
\end{equation}

%%%%%%%%%%%%%%%%%%%%%%%%%%%
\subsubsection{Schubert calculus}
%%%%%%%%%%%%%%%%%%%%%%%%%%%

Let $\iota: U_\C\croc \tU_\C$ be the complexification of $\iota:U\croc \tU$.  To any non-zero rational element $\gamma\in\tgot$, we associate the parabolic subgroups 
\begin{equation}\label{eq:P-gamma}
\tP_\gamma=\{g\in \tU_\C, \lim_{t\to\infty}\exp(-it\gamma)g\exp(it\gamma)\ {\rm exists}\}\quad {\rm and}\quad P_\gamma=\tP_\gamma\cap U_\C.
\end{equation}

We consider the projective varieties $\Fcal_\gamma:=U_\C/P_\gamma$ and $\widetilde{\Fcal}_\gamma:=\tU_\C/\tP_\gamma$, with the canonical embedding 
$\iota:\Fcal_\gamma\croc \widetilde{\Fcal}_\gamma$. Let $B\subset U_\C$ (resp. $\tB\subset \tU_\C$) be the Borel subgroup associated to the choice of the 
Weyl chamber $\tgot_{+}$ (resp. $\ttgot_{+}$).

We associate to $(w,\tw)\in W\times \tW$, the Schubert cells
$$
\tilde{\Xgot}^o_{\tilde{w},\gamma}:= \tB[\tilde{w}]\subset \widetilde{\Fcal}_\gamma\qquad {\rm and}\qquad 
\Xgot^o_{w,\gamma}:= B[w]\subset \Fcal_\gamma.
$$
The corresponding Schubert varieties are $\tilde{\Xgot}_{\tilde{w},\gamma}:=\overline{\tilde{\Xgot}^o_{\tilde{w},\gamma}}$ and $\Xgot_{w,\gamma}:=\overline{\Xgot^o_{w,\gamma}}$. 

We consider the cohomology\footnote{Here, we use singular cohomology with integer coefficients.} rings $H^*(\widetilde{\Fcal}_\gamma,\Z)$ and 
$H^*(\Fcal_\gamma,\Z)$. Let 
$$
\iota^*:H^*(\widetilde{\Fcal}_\gamma,\Z)\to H^*(\Fcal_\gamma,\Z)
$$ 
be the pull-back map in cohomology. If $Y$ is an irreducible closed subvariety of $\widetilde{\Fcal}_\gamma$, we denote by 
$[Y]\in H^{2n_Y}(\widetilde{\Fcal}_\gamma,\Z)$ its cycle class in cohomology : here $n_Y={\rm codim}_\C(Y)$. 
Recall that the cohomology class $[pt]$ associated to a singleton $Y=\{pt\}\subset \Fcal_\gamma$  is a basis of 
$H^{\mathrm{max}}(\Fcal_\gamma,\Z)$.

In the next section we will consider a rational element $\gamma\in\tgot$ and $(w,\tw)\in W\times \tW$ satisfying the relation 
$[\Xgot_{w,\gamma}]\cdot \iota^*([\tilde{\Xgot}_{\tilde{w},\gamma}])= k[pt]$ in $H^*(\Fcal_\gamma,\Z)$, with $k\geq 1$.
This cohomological condition implies in particular that $\dim_\C(\Xgot_{w,\gamma})={\rm codim}_\C(\tilde{\Xgot}_{\tilde{w},\gamma})$ which is equivalent to the relation 
\begin{equation}\label{eq:dimension}
\sharp\left\{\alpha\in\Rgot^+, \langle\alpha,w\gamma\rangle> 0\right\}=\sharp\left\{\tilde{\alpha}\in\widetilde{\Rgot}^+,\langle\tilde{\alpha},\tilde{w}\gamma\rangle< 0\right\}.
\end{equation}

We finish this section by considering the particular case where $U_\C=GL_n(\C)$ is embedded diagonally in  $\tU_\C=GL_n(\C)\times GL_n(\C)$. For $1\leq r<n$, the 
vector $\gamma_r=(\underbrace{-1,\ldots,-1}_{r \ times},0,\ldots,0)\in\R^n\simeq \tgot$ is admissible and the flag manifolds $\Fcal_{\gamma_r}$ and 
$\widetilde{\Fcal}_{\gamma_r}$ admits a canonical identifications respectively with the Grassmanians $\G(r,n)$ and $\G(r,n)\times\G(r,n)$. The map 
$\iota:\Fcal_\gamma\croc \widetilde{\Fcal}_\gamma$ corresponds to the diagonal embedding $\iota:\G(r,n)\to\G(r,n)\times\G(r,n)$.

For $w\in W\simeq \Sgot_n$, the Schubert variety $\Xgot_{w,\gamma_r}\subset \G(r,n)$, which depends only of the subset $K=w([r])\subset [n]$, is denoted $\Xgot_K$. 
Similarly, for $w=(w_1,w_2)\in \widetilde{W}\simeq \Sgot_n\times\Sgot_n$, the Schubert variety $\tilde{\Xgot}_{\tilde{w},\gamma_r}\subset \G(r,n)\times \G(r,n)$ is equal to 
$\Xgot_I\times\Xgot_J$, where $I=w_1([r])$ and $J=w_2([r])$.

In this setting, we have the following classical result.

\begin{lemma}
The following statements are equivalent:
\begin{itemize}
\item $[\Xgot_{w,\gamma}]\cdot \iota^*([\tilde{\Xgot}_{\tilde{w},\gamma}])= \ell[pt]$ in $H^*(\Fcal_\gamma,\Z)$, with $\ell\geq 1$.
\item $[\Xgot_I]\cdot[\Xgot_J]\cdot[\Xgot_K]=\ell[pt]$ in $H^*(\G(r,n),\Z)$, with $\ell\geq 1$.
\item The Littlewood-Richardson coefficient $c^{\mu(K)}_{\mu(I^o),\mu(J^o)}$ is non-zero.
\end{itemize}
\end{lemma}

%%%%%%%%%%%%%%%%%%%%%%%%%%%%%%%%%%%%
\subsubsection{Description of $\LR(U,\tU)$}\label{sec:LR-U-tildeU}
%%%%%%%%%%%%%%%%%%%%%%%%%%%%%%%%%%%%

We can finally describe the cone $\LR(U,\tU)$.

\begin{theorem}\label{theo:LR-general} Let $(\xi,\txi)\in\tgot_{+}\times \ttgot_{+}$. We have  $U\xi\subset \pi\big(\,\tU\txi\,\big)$ if and only if 
\begin{equation}\label{eq:inequality-polytope}
\langle \tilde{\xi},\tilde{w}\gamma\rangle\geq \langle \xi,w_ow\gamma\rangle
\end{equation}
for any $(\gamma,w,\tilde{w})\in\tgot \times W\times \tilde{W}$ satisfying the following properties:
\begin{enumerate}
\item[a)] $\gamma$ is admissible antidominant.
\item[b)] $[\Xgot_{w,\gamma}]\cdot \iota^*([\tilde{\Xgot}_{\tilde{w},\gamma}])= [pt]$ in $H^*(\Fcal_\gamma,\Z)$.
\item[c)] Identity (\ref{eq:trace-condition-w-tilde}) holds.
\end{enumerate}
The result still holds if we replace {\rm b)} by the weaker condition
$$
{\rm b')} \qquad [\Xgot_{\gamma}]\cdot \iota^*([\tilde{\Xgot}_{\tilde{w},\gamma}])= \ell[pt], \ell\geq 1\qquad {\rm in}\quad H^*(\Fcal_\gamma,\Z).
$$
\end{theorem}

\begin{remark}\label{rem:condition-c}
Suppose that there exists $c_\gamma>0$ such that 
$|\langle\alpha,w\gamma\rangle|\,$ and $\, | \langle\tilde{\alpha},\tilde{w}\gamma\rangle|\,$ belongs to $\{0, c_\gamma\}$, $\forall (w,\tw)\in W\times \tilde{W}$, 
$\forall (\alpha,\tilde{\alpha})\in \Rgot(\ugot)\times\Rgot(\tugot)$. Then condition {\rm c)} follows from condition {\rm b)} (see (\ref{eq:dimension})). 
\end{remark}

When the closed connected subgroups $\iota: U\croc \tU$ satisfy Assumption \ref{hypothese} the subspace $\zgot:=\{X\in \tgot, \alpha(X)=0, \forall \alpha\in\Sigma(\tugot/\ugot)\}$ is equal to $Z_{\tugot}\cap\ugot$. Let $\tgot=Z_{\tugot}\cap\ugot\oplus\tgot_1$ be a rational decomposition. 
Any rational element $\gamma\in\tgot$ can be written $\gamma=\gamma_0+\gamma_1$ where $\gamma_0\in Z_{\tugot}\cap\ugot$ and $\gamma_1\in\tgot_1$ are rational. 
Two cases occur :
\begin{itemize}
\item If $\gamma_1=0$, then $\gamma$ satisfies conditions a), b) and c). The inequalities (\ref{eq:inequality-polytope}) given by these central elements shows that $\tilde{\xi}-\xi$ is 
orthogonal to  $Z_{\tugot}\cap\ugot$.

\item If $\gamma_1\neq 0$ then it is immediate to see that $\gamma$ satisfies a), b) and c) if and only if $\gamma_1$ does also.  Moreover, as 
$\tilde{\xi}-\xi$ is orthogonal to  $Z_{\tugot}\cap\ugot$, $\langle \tilde{\xi},\tilde{w}\gamma\rangle\geq \langle \xi,w_ow\gamma\rangle$ if and only if 
$\langle \tilde{\xi},\tilde{w}\gamma_1\rangle\geq \langle \xi,w_ow\gamma_1\rangle$.
\end{itemize}

Thus, when Assumption \ref{hypothese} is satisfied, $\LR(U,\tU)$ is described by the condition $\tilde{\xi}-\xi\in(Z_{\tugot}\cap\ugot)^\perp$ and a finite number of inequalities of the form 
$\langle \tilde{\xi},\tilde{w}\gamma_1\rangle\geq \langle \xi,w_ow\gamma_1\rangle$.

\medskip

Many people have contributed to Theorem \ref{theo:LR-general}. The first input was given by Klyachko \cite{Klyachko} with a refinement by 
Belkale \cite{Belkale06}, in the case of $SL_n\croc (SL_n)^s$. The case $U_\C\croc (U_\C)^s$ has been treated by Belkale-Kumar \cite{BK06} and by Kapovich-Leeb-Millson 
\cite{KLM-memoir-08}. Recall that {\em Condition c)} is related to the notion of Levi-movability introduced by Belkale-Kumar \cite{BK06}.  
Finally, Berenstein-Sjamaar \cite{Berenstein-Sjamaar} and Ressayre \cite{Ressayre10, Ressayre11}  have studied the general case. Ressayre \cite{Ressayre10} also proved the irredundancy of the list of inequalities.

We refer the reader to the survey articles \cite{Fulton00,Brion-Bourbaki,Kumar14}  for details.

%%%%%%%%%%%%%%%%%%%%%%%%%%%%%
\subsection{The cone $\LR(m,n)$} \label{sec:LR-m-n}
%%%%%%%%%%%%%%%%%%%%%%%%%%%%%

Let $m,n\geq 1$. Let us write an Hermitian matrix $X\in Herm(m+n)$ by blocks 
$X=\begin{pmatrix}
X_{\bf I}& *\\
*& X_{\bf II}
\end{pmatrix}$
where $X_{\bf I}\in Herm(m)$ and $X_{\bf II}\in Herm(n)$. In this section, we are interested in the cone $\LR(m,n):=\left\{\left(\e(X),\e(X_{\bf I}),\e(X_{\bf II})\right); \ X\in Herm(m+n)\right\}$. 
Thanks to Theorem \ref{theo:LR-general}, we obtain the following description of $\LR(m,n)$. The details of the proof are given in the next section.

\begin{theorem}\label{theo:LR-m-n} The triplet $(x,y,z)\in \R^{m+n}_+\times \R^{m}_+\times \R^{n}_+$ belongs to $\LR(m,n)$ if and only if the following conditions hold:
\begin{itemize}
\item $|x|=|y|+|z|$,
\item $x_{n+k}\leq y_k\leq x_k$, $\forall k\in [m]$,
\item $x_{m+\ell}\leq z_\ell\leq x_\ell$, $\forall \ell\in [n]$,
\item $|x|_A\geq |y|_B + |z|_C$, for any triplet $A,B,C$ satisfying:
\begin{enumerate}
\item $B\subset [m]$ and $C\subset [n]$ are strict subsets,
\item $A\subset [m+n]$ and $\sharp A =\sharp B +\sharp C$,
\item the Littlewood-Richardson coefficient $c^{\mu(A)}_{\mu(B),\mu(C)}$ is non-zero.
\end{enumerate}
\end{itemize} 
Moreover, the condition $c^{\mu(A)}_{\mu(B),\mu(C)}\neq 0$ is equivalent to $(\mu(A),\mu(B),\mu(C))\in \LR(u,v)$, where $u=\sharp B$ and $v=\sharp C$.
\end{theorem}

\begin{remark}
In Theorem \ref{theo:LR-m-n}, we can strenghten condition 3. by requiring that \break $c^{\mu(A)}_{\mu(B),\mu(C)}=1$.
\end{remark}

\begin{remark}
In \cite{Li-Poon03}, Li and Poon also obtained a characterization of the cone $\LR(m,n)$ by means of the following inequalities: 
$|x|_I\leq |y|_{J\cap[m]} + |z|_{K\cap[n]}$,\ $\forall (I,J,K)\in \LR^{n+m}_r$,\  $\forall r<n+m$.
\end{remark}

We will see in the next section that $c^{\mu(A)}_{\mu(B),\mu(C)}\neq 0$ if and only $c^{\mu((A^o)^c)}_{\mu((B^o)^c)),\mu((C^o)^c))}\neq 0$. Since the relation 
$|x|_{(A^o)^c}\geq |y|_{(C^o)^c} + |z|_{(C^o)^c}$ is equivalent to $|x|_{A^o}\leq |y|_{C^o} + |z|_{C^o}$, in Theorem \ref{theo:LR-m-n}, we can rewrite  the last condition by requiring that 
\begin{equation}\label{eq:duality-equations}
|x|_A\geq |y|_B + |z|_C\quad {\rm and}\quad |x|_{A^o}\leq |y|_{C^o} + |z|_{C^o}
\end{equation}
for all strict subsets $A\subset [m+n]$, $B\subset [m]$, $C\subset [n]$ that satisfy $\sharp A =\sharp B +\sharp C\leq \frac{1}{2}(m+n)$ and $c^{\mu(A)}_{\mu(B),\mu(C)}\neq 0$.

%%%%%%%%%%%%%%%%%%%%%%%%%%%%%
\subsection{Proof of Theorem \ref{theo:LR-m-n}}
%%%%%%%%%%%%%%%%%%%%%%%%%%%%%

We work with the unitary group $\tU=U_{m+n}$ and the subgroup $U=U_m\times U_n$ embedded diagonally. We consider the orthogonal projection
$\pi_0:Herm(m+n)\to Herm(m)\times Herm(n)$ that sends $X$ to $\pi_0(X)=(X_{\bf I},X_{\bf II})$. The cone $\LR(m,n)$ is formed by the triplets 
$(x,y,z)\in \R^{m+n}_+\times \R^{m}_+\times \R^{n}_+$ satisfying 
$$
U_m\cdot\diag(y)\times U_n\cdot\diag(z)\subset\pi_0\left(U_{m+n}\cdot\diag(x)\right).
$$  
Thus $\LR(m,n)=\LR(U_m\times U_n,U_{m+n})$.
%%%%%%%%%%%%%%%%%%%
\subsubsection{Admissible elements}
%%%%%%%%%%%%%%%%%%%

We work with the maximal torus $T\subset U$ of diagonal matrices. The set of roots relatively to the action of $T$ on $\tugot/\ugot\simeq M_{m,n}(\C)$ is 
$\Sigma:=\{e^*_i-f^*_j; i\in [m], j\in [n]\}$. 

The center of $\tugot$ is generated by $\gamma_o:=(1,\ldots,1)\in\R^{m+n}\simeq \tgot$. For any $(r,s)\in \{0,\ldots,m\}\times \{0,\ldots,n\}$, we define 
$$
\gamma_{r,s}=(\underbrace{-1,\ldots,-1}_{r \ times},0,\ldots,0)\oplus (\underbrace{-1,\ldots,-1}_{s\  times},0,\ldots,0)\in \R^{m}\times\R^n\simeq \tgot.
$$

\begin{lemma} Let $\gamma\in \tgot$ be an admissible element. There exists $(a,b)\in \Q\times \Q^{\geq 0}$, $(w,w')\in\Sgot_m\times\Sgot_n$, and $(r,s)$ such that 
$\gamma=a\gamma_o+ b (w,w')\gamma_{r,s}$. The couple $(r,s)$ must satisfy the auxiliary conditions: either $0\!<\!r\!<\!m$ and $0\!<\!s\!<\!n$ or 
$(r,s)\in\{(1,0),(0,1),(m-1,n),(n,m-1)\}$.
\end{lemma}

{\em Proof :}   Consider an admissible vector $\gamma=(\gamma_1,\ldots,\gamma_m;\gamma'_1,\ldots,\gamma'_n)$ that is linearly independent to $\gamma_o$. The relation 
${\rm Vect}(\Sigma\cap\gamma^\perp)={\rm Vect}(\Sigma)\cap\gamma^\perp$ means that $(\Sigma\cap\gamma^\perp)^\perp$ is a subspace of dimension $2$.
Here $\Sigma\cap\gamma^\perp$ is the set of vectors $e^*_i-f^*_j$ such that $\gamma_i=\gamma_j'$. For $\alpha\in\R$, we define $[m]_\alpha:=\{i\in[m], \gamma_i=\alpha\}$ and 
$[n]_\alpha:=\{j\in[n], \gamma'_j=\alpha\}$. Hence $\Sigma\cap\gamma^\perp$ is parameterized by $\coprod_{\alpha\in L} [m]_\alpha\times [n]_\alpha$ where 
$L=\{\alpha\in \R, [m]_\alpha\neq \emptyset \ {\rm and} \ [n]_\alpha\neq \emptyset \}$ is a finite set. 

Consider first the case where $\cup_{\alpha\in L}[m]_\alpha\neq [m]$. Let $k\notin \cup_{\alpha\in L}[m]_\alpha$. Then 
$(\Sigma\cap\gamma^\perp)^\perp$, which is of dimension $2$, contains the vectors $\gamma_o, \gamma$ and $e_k$. Hence, 
$\gamma$ is a linear combinaison of $\gamma_o$ and $e_k$: we check easily that there exists $(a,b)\in \Q\times \Q^{> 0}$, $w\in\Sgot_m$, and $(r,s)\in\{(1,0),(m-1,n)\}$ such that 
$\gamma=a\gamma_o+ b w\gamma_{r,s}$. 

If $\cup_{\alpha\in L}[n]_\alpha\neq [n]$, we prove similarly that $\gamma=a\gamma_o+ b w'\gamma_{r,s}$ for $(a,b)\in \Q\times \Q^{> 0}$, $w'\in\Sgot_n$, and 
$(r,s)\in\{(0,1),(m,n-1)\}$. 

Let us consider the last case where $\cup_{\alpha\in L}[m]_\alpha= [m]$ and $\cup_{\alpha\in L}[n]_\alpha= [n]$. Then $\gamma=\sum_{\alpha\in L} \alpha V_\alpha$ with  
$V_\alpha=\sum_{i\in [m]_\alpha, j\in [n]_\alpha} e_i+f_j$. The vectors $\{V_\alpha,\alpha\in L\}$ define an independent family of the subspace $
(\Sigma\cap\gamma^\perp)^\perp$ which is of dimension $2$, so $\sharp L\leq 2$. Since  $\gamma$ is linearly independent to $\gamma_o$, 
the set $L$ has cardinal $2$. Now we see that there exists $(a,b)\in \Q\times \Q^{> 0}$, $(w,w')\in\Sgot_m\times\Sgot_n$, and $r<m, s<n$ such that 
$\gamma=a\gamma_o+ b (w,w')\gamma_{r,s}$. $\Box$

\medskip

Here, the remark \ref{rem:condition-c} applies, so  condition c) will follow from condition b).

%%%%%%%%%%%%%%%%%%%%%%%%%%%%%%%
\subsubsection{Cohomological conditions and inequalities}
%%%%%%%%%%%%%%%%%%%%%%%%%%%%%%%

The Lie algebra $\tgot$ is identified with $\R^{m+n}\simeq \R^m\times\R^n$. 

{\bf First case:} The two vectors $\pm \gamma_o$ are admissible elements, and satisfy conditions $a)$, $b)$ and $c)$ 
of Theorem \ref{theo:LR-general} in an obvious way. In this cases, the corresponding inequalities 
$\pm(x,\gamma_o)\geq \pm((y,z),w_o\gamma_o)$ are equivalent to $|x|=|y|+|z|$.

{\bf Second case:} We work now with the admissible element $\gamma_{r,s}$ in the situation where $r\in [m-1]$ and $s\in [n-1]$. The flag manifold $GL_m(\C)\times GL_n(\C)/P_{\gamma_{r,s}}$ admits a natural identification with the product of Grassmannians $\G(r,m)\times \G(s,n)$. Similarly, the flag manifold $GL_{m+n}(\C)/\widetilde{P}_{\gamma_{r,s}}$ is isomorphic to the Grassmannian $\G(r+s,m+n)$. 
The map $\iota_\C : GL_m(\C)\times GL_n(\C) \to GL_{m+n}(\C)$  factorises to a smooth map 
 $\iota_{r,s} :  \G(r,m)\times\G(s,n)\to  \G(r+s,m+n)$ defined by $\iota_{r,s}(V_1,V_2)= V_1\oplus V_2$. 
 
Let $w=(w_1,w_2)\in W\simeq \Sgot_m\times\Sgot_n$ and let $B=w_1([r])\subset [m]$ and $C=w_2([s])\subset [n]$ be the corresponding subsets. The associated Schubert variety is 
 $\Xgot_{w,\gamma_{r,s}}= \Xgot_B\times \Xgot_C\subset  \G(r,m)\times\G(s,n)$. 
 
 In the same way, to $\tilde{w}\in\tilde{W}\simeq \Sgot_{m+n}$, we associate the subset $A=\tilde{w}([r])\cup \tilde{w}([s]+m)\subset [m+n]$ and the Schubert variety 
 $\tilde{\Xgot}_{\tilde{w},\gamma_{r,s}}= \Xgot_A \subset  \G(r+s,m+n)$.
 
 \begin{lemma} The following identities are equivalent:
 \begin{enumerate}
\item $[\Xgot_{w,\gamma}]\cdot \iota^*([\tilde{\Xgot}_{\tilde{w},\gamma}])=\ell [pt], \ell\geq 1$,
\item $c^{\mu(A^o)}_{\mu(B),\mu(C)}=\ell\geq 1$,
\item $c^{\mu(A^c)}_{\mu((B^o)^c),\mu((C^o)^c)}= \ell\geq 1$.
\end{enumerate}
 \end{lemma}

{\em Proof:} Recall that we associate a partition $\lambda(A)=(\lambda_1\geq\cdots\lambda_r)$ with a subset $A=\{a_1<i_2<\cdots<a_r\}\subset [n]$ 
of cardinality $r$, by posing  $\lambda_k=n-r+k-a_k$, $\forall k\in [r]$.

Let $\bigwedge_r[x]=\Z[x_1,\ldots,x_r]^{\Sgot_r}$  be the ring of symmetric polynomials, with integral coefficients, in $r$ variables. 
For any partition $\nu$ of length $r$, we associate its Schur polynomial ${\bf s}_\nu(x)\in\bigwedge_r[x]$. The family $({\bf s}_\nu)$ determine a $\Z$-basis of $\bigwedge_r[x]$.

Let us recall recall the following classical fact (see \S 3.2.2 in \cite{Manivel-98}). The map $\phi_{r}:\bigwedge_r[x]\longrightarrow H^{*}(\G(r,m))$ defined by the relations 
$$
\phi_{r}({\bf s}_\nu)=
\begin{cases}
\sigma_\nu\hspace{8mm} {\rm if}\quad \nu_1\leq m-r,\\
0\hspace{1cm}{\rm if}\quad \nu_1> m-r.
\end{cases}
$$
is a ring morphism. Here $\sigma_\nu$ denotes the cohomology class $[\Xgot_{D}]$ defined by a subset $D\subset [m]$ of cardinality $r$ such that $\nu=\lambda(D)$. 
In the same way we consider the ring  $\bigwedge_{r+s}[x,y]=\Z[x_1,\ldots,x_r,y_1,\ldots,y_s ]^{\Sgot_{r+s}}$ and the morphism
$\phi_{r+s}:\bigwedge_{r+s}[x,y]\longrightarrow H^{*}(\G(r+s,m+n))$.  Let us denote by $R: \bigwedge_{r+s}[x,y]\to $
$ \bigwedge_r[x] \otimes  \bigwedge_s[y]$ the restriction morphism. It is not hard to check that the following diagram is commutative:
$$
\begin{tikzcd}
  \bigwedge_r[x] \arrow[d, "\phi_r"] &\otimes & \bigwedge_s[y] \arrow[d, "\phi_s"]  & \bigwedge_{r+s}[x,y]\arrow[l, "R"]\arrow[d, "\phi_{r+s}"] \\
  H^{*}(\G(r,m)) &\otimes &  H^{*}(\G(s,n)) &H^{*}(\G(r+s,m+n))\arrow[l, "j^*"].
\end{tikzcd}
$$

As $[\Xgot_{w,r}]=\sigma_{\lambda(B)}\otimes\sigma_{\lambda(C)}$ and $[\tilde{\Xgot}_{\tw,r}]=\sigma_{\lambda(A)}$, the previous diagram tell us that the integer $\ell$ such that 
$[\Xgot_{w,r}]\cdot \iota_r^*([\tilde{\Xgot}_{\tw,r}])= \ell [pt]$ is equal to the coefficient of $R({\bf s}_{\lambda(A)}(x,y))$ relatively to 
${\bf s}_{\lambda(B^o)}(x)\otimes {\bf s}_{\lambda(C^o)}(y)$ : in other words $\ell$ is equal to the Littlewood-Richardson coefficient  
$c^{\lambda(A)}_{\lambda(B^o),\lambda(C^o)}=c^{\mu(A^o)}_{\mu(B),\mu(C)}$ (see \cite{MacDonald}, \S I.5). The equivalence between {\em 1.} and {\em 2.} is proved. 

\medskip

Let us consider $r', s'$ such that $r+r'=m$ and $s+s'=n$. The canonical bilinear form on $\C^{m+n}$ permits to define the map $\delta: \G(r'+s',m+n)\to \G(r+s,m+n)$ 
that sends a subspace $F\subset\C^{m+n}$ to its orthogonal $F^\perp$. Let $\delta^*: H^{*}(\G(r+s,m+n))\to H^{*}(\G(r'+s',m+n))$ denote the pullback map in cohomology. If we consider similar maps $\delta^*: H^{*}(\G(r,m))\to H^{*}(\G(r',m))$ and $\delta^*: H^{*}(\G(s,n))\to H^{*}(\G(s',n))$, we have a commutative diagram:
$$
\begin{tikzcd}
  H^{*}(\G(r,m)) \arrow[d, "\delta^*"] &\otimes & H^{*}(\G(s,n)) \arrow[d, "\delta^*"]  & H^{*}(\G(r+s,m+n))\arrow[l, "j^*"]\arrow[d, "\delta^*"] \\
  H^{*}(\G(r',m)) &\otimes &  H^{*}(\G(s',n)) &H^{*}(\G(r'+s',m+n))\arrow[l, "j^*"].
\end{tikzcd}
$$

This allows us to see that $\sigma_{\lambda(B)}\otimes\sigma_{\lambda(C)}\cdot j^*(\sigma_{\lambda(A)})= \ell[pt],k\geq 1$ if and only if 
$\delta^*(\sigma_{\lambda(B)})\otimes\delta^*(\sigma_{\lambda(C)})\cdot j^*(\delta^*(\sigma_{\lambda(A)}))= \ell[pt],k\geq 1$. Since we have
$\delta^*(\sigma_{\lambda(X)})=\sigma_{\mu(X^c)}$ as a general rule, the previous relation is equivalent to $c^{\mu(A^c)}_{\mu((B^o)^c),\mu((C^o)^c)}= \ell\geq 1$. 
The equivalence between {\em 2.} and {\em 3.} is proved.  $\Box$

\medskip

The inequalities associated to $\gamma_{r,s}$ are 
$$
-|x|_A=\langle x,\tilde{w}\gamma_{r,s}\rangle\geq \langle (y,z),w_ow\gamma\rangle=-|y|_{B^o} -|z|_{C^o} 
$$

Using $|x|=|y|+|z|$, we obtain $|x|_{A^c}\geq |y|_{(B^o)^c} +|z|_{(C^o)^c}$ for any strict subsets $A\subset [m+n]$, $B\subset [m]$ and $C\subset [n]$ satisfying 
$\sharp A =\sharp B +\sharp C$ and $c^{\mu(A^c)}_{\mu((B^o)^c),\mu((C^o)^c)}= \ell\geq 1$.

\medskip

{\bf Third case: $(r,s)\in\{(1,0),(0,1),(m-1,n),(n,m-1)\}$.} Here we use the same type of argument as before.

$(r,s)=(1,0)$: we obtain the inequalities $x_{n+k}\leq y_k,\forall k\in[m]$.

$(r,s)=(0,1)$: we obtain the inequalities $x_{m+\ell}\leq z_\ell,\forall \ell\in[n]$.

$(r,s)=(m-1,n)$: we obtain the inequalities $y_{k}\leq x_k,\forall k\in[m]$.

$(r,s)=(m,n-1)$: we obtain the inequalities $z_{\ell}\leq x_\ell,\forall \ell\in[n]$.

\medskip
 
The proof of Theorem \ref{theo:LR-m-n} is completed. $\Box$

%%%%%%%%%%%%%%%%%%%%%%%%%%%%%%%%
\subsection{A consequence of the O'Shea-Sjamaar Theorem}
%%%%%%%%%%%%%%%%%%%%%%%%%%%%%%%%

%%%%%%%%%%%%%%%%%%%%%%%%%%%%%%%%
\subsubsection{First setting: compact Lie groups with involution}
%%%%%%%%%%%%%%%%%%%%%%%%%%%%%%%%

Let $\tU$ be a compact connected Lie group equipped with an involution $\sigma$. The Lie algebra of $\tU$ admit the decomposition 
$\tugot=\tugot^{\sigma}\oplus\tugot^{-\sigma}$ that is invariant under the action of the subgroup $\tU^\sigma$. 
We start with a basic but important fact (see \cite{OSS}, Example 2.9).

\begin{lemma}\label{lemma:orbit-intersection}
For any adjoint orbit $\tOcal\subset\tugot$, the intersection $\tOcal\cap\tugot^{-\sigma}$ is either empty or an orbit of the connected subgroup $\tK:=(\tU^\sigma)_0$.
\end{lemma}

Let $U\subset\tU$ be a subgroup invariant under $\sigma$. Let us choose an invariant scalar product $(-,-)$ on the Lie algebra $\tugot$ of $\tU$ such that $\sigma\in O(\tugot)$.
At the level of Lie algebras, we consider the orthogonal projection $\pi:\tugot\to\ugot$ relatively to the scalar product $(-,-)$.

One of the main tool used in this paper is the following result, which is a consequence of the O'Shea-Sjamaar Theorem (see \cite{OSS}, Section 3). Let $K$ be the connected component of 
$U^\sigma$.

\begin{proposition}\label{prop:oshea-sjamaar-1}
Let $\xi\in\ugot^{-\sigma}$ and $\txi\in\tugot^{-\sigma}$. The following conditions are equivalent:
\begin{enumerate}
\item $U\xi\subset \pi\left(\tU\txi\right)$,
\item $K\xi\subset \pi\left(\tK\txi\right)$.
\end{enumerate}
\end{proposition}

%%%%%%%%%%%%%%%%%%%%%%%%%%%%%%%%
\subsubsection{Second setting: real reductive Lie groups}
%%%%%%%%%%%%%%%%%%%%%%%%%%%%%%%%

Let $\iota:G\croc \tG\subset GL_N(\R)$ be two connected real reductive Lie groups admitting a complexification $\iota_\C:G_\C\croc \tG_\C\subset GL_{N}(\C)$. 
It is for example the case when $G$ and $\tG$ are semisimple (see \cite{Knapp-book}, \S VII.1). Let us denote by
\begin{itemize}
\item $K=G\cap SO_N(\R)$ and $\tK=\tG\cap SO_N(\R)$ the maximal compact subgroups of $G$ and $\tG$. Their Lie algebras are denoted by 
$\iota:\kgot\croc\tkgot$.
\item $U= G_\C\cap U_{N}$ and  $\tU=\tG_\C\cap U_N$ the maximal compact subgroups of $G_\C$ and $\tG_\C$. Their Lie algebras are denoted by 
$\iota:\ugot\croc\tugot$.
\end{itemize}

Consider the Cartan decompositions, $\ggot=\kgot\oplus \pgot$ and $\tggot=\tkgot\oplus \tpgot$, of $G$ and $\tG$. At the level of Lie algebras, we have $\tugot=\tkgot\oplus i \tpgot$ and $\ugot=\kgot\oplus i \pgot$. The antilinear conjugation on $GL_N(\C)$ defines an involution 
$\sigma$ on $U\croc \tU$ such that $K\croc\tK$ are respectively equal to the connected components of $U^\sigma$ and $\tU^\sigma$. We see also that 
$\ugot^{-\sigma}=i\pgot$ and $\tugot^{-\sigma}=i\tpgot$.

Let $\pi:\tggot_\C\to\ggot_\C$ be the orthogonal projection relatively to the Hermitian norm ${\rm Tr}(X^*X)^{1/2}$ on $\glgot_N(\C)$. 

\begin{proposition}\label{prop:oshea-sjamaar-2}
Let $X\in\pgot$ and $\tX\in\tpgot$. The following conditions are equivalent:
\begin{enumerate}
\item $UX\subset \pi\left(\tU\tX\right)$,
\item $KX\subset \pi\left(\tK\tX\right)$.
\end{enumerate}
\end{proposition}

{\em Proof:} It follows from Proposition \ref{prop:oshea-sjamaar-1} and the fact that the projection $\pi:\tggot_\C\to\ggot_\C$ is complex linear. $\Box$

%%%%%%%%%%%%%%%%%%%%%%%%%%%%%%%%%%%%%%%%%%%%%%%%
%%%%%%%%%%%%%%%%%%%%%%%%%%%%%%%%%%%%%%%%%%%%%%%%
\section{The cone $\Acal(p,q)$}\label{sec:A-p-q}
%%%%%%%%%%%%%%%%%%%%%%%%%%%%%%%%%%%%%%%%%%%%%%%%
%%%%%%%%%%%%%%%%%%%%%%%%%%%%%%%%%%%%%%%%%%%%%%%%

Here, we work with the reductive real Lie group $U(p,q)=\{g\in GL_n(\C),\, g^*I_{p,q}g=I_{p,q}\}$,  where $I_{p,q}=\diag(I_p, -I_q)$.

%%%%%%%%%%%%%%%%%%%%%%%%%%%%%%%
\subsection{Matrix identities}
%%%%%%%%%%%%%%%%%%%%%%%%%%%%%%%

Let us decompose a $n$-square hermitian matrix 
$X=\begin{pmatrix}
X_{11}& X_{12}\\
X_{12}^*& X_{22}
\end{pmatrix}$
by blocks, where $X_{12}\in M_{p,q}(\C)$. Recall that $(\lambda, s)\in \Acal(p,q)$ if and only if there exists an hermitian matrix $X$ such that $\lambda=\e(X)$ and $s=\s(X_{12})$. Let us consider 
$$
\widetilde{X}=-I_{p,q}X I_{p,q}=\begin{pmatrix}
-X_{11}& X_{12}\\
X_{12}^*& -X_{22}
\end{pmatrix}
\quad 
{\rm and}
\quad 
X+\widetilde{X}=2\begin{pmatrix}
0& X_{12}\\
X_{12}^*& 0
\end{pmatrix}.
$$
If we look at the eigenvalues of this three Hermitian matrices, we obtain, following the notations of Section \ref{sec:augmented}, 
$$
\e(\widetilde{X})=\e(X)^*\quad  {\rm and} \quad \e(X+\widetilde{X})=2\,\widehat{\s(X_{12})}^{p,q},
$$
i.e. $(\e(X),\e(X)^*,2\,\widehat{\s(X_{12})}^{p,q})\in \horn(n)$.

From the above identities, we see that 
any $(\lambda, s)\in \Acal(p,q)$ satisfies the relation
\begin{equation}\label{eq:OSS-relation-A-p-q}
(\lambda,\lambda^*, 2\,\widehat{s}\,^{p,q})\in \horn(n).
\end{equation}

In the following sections, we explain how the O'Shea-Sjamaar theorem (see Proposition \ref{prop:oshea-sjamaar-2}) allows us to see that relation (\ref{eq:OSS-relation-A-p-q}) characterizes the cone $\Acal(p,q)$.  

%%%%%%%%%%%%%%%%%%%%%%%%%%%%%%%
\subsection{Complexification and antiholomorphic involution}
%%%%%%%%%%%%%%%%%%%%%%%%%%%%%%%

We work with the reductive real Lie groups $G:= U(p,q)$ and $\tG:= GL_n(\C)$. Let us denote by $\iota: G\to \tG$ the canonical embedding. The unitary group 
$\tK:=U_n$ is a maximal compact subgroup of $\tG$. Let $\tpgot:= Herm(n)\subset \glgot_n(\C)$ be the subspace of Hermitian matrices.

The subgroup $K:=\tK\cap U(p,q)\simeq U_p\times U_q$ is a maximal compact sugroup of $G$, and the map 
$Y\mapsto \widehat{Y}^{p,q}=\begin{pmatrix} 0_{pp}& Y\\ Y^*& 0_{qq}\end{pmatrix}$ defines an identification between $M_{p,q}(\C)$ and 
the subspace $\pgot:=\tpgot\cap\ggot$.

The complexification of the group $G$ is $G_\C:= GL_n(\C)$. We consider the antiholomorphic involution $\sigma$ on $G_\C$ defined by 
$\sigma(g)= I_{p,q}(g^*)^{-1}I_{p,q}$. The subgroup $G$ is the fixed point set of $\sigma$.

The complexification of the group $\tG$ is $\tG_\C:= GL_n(\C)\times GL_n(\C)$. The inclusion $\tG\croc \tG_\C$ is given by the map $g\mapsto (g,\bar{g})$. 
We consider the antiholomorphic involution $\widetilde{\sigma}$ on $\tG_\C$ defined by $\widetilde{\sigma}(g_1,g_2)= (\overline{g_2},\overline{g_1})$. 
The subgroup $\tG$ corresponds to the fixed point set of $\widetilde{\sigma}$. The embedding $\iota: G\croc \tG$ admits a complexification 
$\iota_\C: G_\C\croc \tG_\C$ defined by $\iota_\C(g)=(g,\overline{\sigma(g)})$: notice that 
$\iota_\C\circ \sigma= \widetilde{\sigma}\circ \iota_\C$.

The groups $U=U_n$ and $\tU= U_n\times U_n$  are respectively  maximal compact sugroups of $G_\C$ and $\tG_\C$. The embedding $\iota_\C: U\croc \tU$ is defined by $\iota_\C(k)=(k,I_{p,q}\bar{k}I_{p,q})$. The fixed point subgroups of the involutions are $U^\sigma=K$ and $\tU^{\tilde{\sigma}}=\tK$.

At the level of Lie algebra, we have a morphism $\iota_\C: \glgot_n(\C)\croc \glgot_n(\C)\times \glgot_n(\C)$ defined by $\iota_\C(X)=(X,\overline{\sigma(X)})$, where 
$\sigma(X)=-I_{p,q}X^*I_{p,q}$.

%%%%%%%%%%%%%%%%%%%%%%%%%%%%%%%
\subsection{Orthogonal projection of orbits}
%%%%%%%%%%%%%%%%%%%%%%%%%%%%%%%

We use on $\glgot_n(\C)\times \glgot_n(\C)$ the euclidean norm $\|(X,Y)\|^2={\rm Tr}(XX^*)+{\rm Tr}(YY^*)$. The subspace orthogonal to the image of $\iota_\C$ is 
$\{(X,-\overline{\sigma(X)}), X\in\glgot_n(\C)\}$. Hence the orthogonal projection
$$
\pi : \glgot_n(\C)\times \glgot_n(\C)\longrightarrow \glgot_n(\C),
$$
is defined by the relations $\pi(X,Y)=\tfrac{1}{2}(X+\overline{\sigma(Y)})$. Note that $\pi$ commutes with the involutions : $\pi\circ \widetilde{\sigma}=\sigma\circ \pi$.

If $X\in Herm(n)$, the corresponding adjoint orbit $U_n\cdot X$, which is entirely determined by the spectrum $\e(X)$, is denoted by $\Ocal_{\e(X)}$. 
If $\lambda=(\lambda_1,\cdots,\lambda_n)$, we denote by $\lambda^*$ the vector $(-\lambda_n,\cdots,-\lambda_1)$: we see that $\e(-X)=\e(X)^*$ 
for any $X\in Herm(n)$. 

The subspace $\tpgot\subset\tggot$ is identified with $\{(X,\overline{X}),X\in Herm(n)\}\subset \glgot_n(\C)\times \glgot_n(\C)$. For any $X\in Herm(n)$, 
the image by the projection $\pi$ of the orbit $\tU\cdot (X,\overline{X})$ is equal to 
$$
\tfrac{1}{2}\left(U_n\cdot X+U_n\cdot \sigma(X)\right)=\tfrac{1}{2}\left(U_n\cdot X+U_n\cdot (-X)\right)=\tfrac{1}{2}\left(\Ocal_{\lambda}+\Ocal_{\lambda^*}\right),
$$
where $\lambda=\e(X)$.

If $Y\in M_{p,q}(\C)$ has singular spectrum $s\in\R^q_{++}$, the spectrum of the Hermitian matrix $\widehat{Y}^{p,q}$ is equal to $\widehat{s}\,^{p,q}$, hence $U\cdot \widehat{Y}^{p,q}$ is equal to $\Ocal_{\widehat{s}\,^{p,q}}$. At this stage, we have proved that for any $(X,Y)\in Herm(n)\times M_{p,q}(\C)$ the following statements are equivalents:
\begin{itemize}
\item $U\cdot\widehat{Y}^{p,q}\subset\pi\left(\tU\cdot (X,\overline{X})\right)$,
\item $2\Ocal_{\widehat{s}\,^{p,q}}\subset \Ocal_{\lambda}+\Ocal_{\lambda^*}$,
\item $(\lambda,\lambda^*, 2\widehat{s}\,^{p,q})\in\horn(n)$,
\end{itemize}
where $\lambda=\e(X)$ and $s=\s(Y)$.

The group $K\simeq U_p\times U_q$ acts canonically $M_{p,q}(\C)\simeq \pgot$. For any $Y\in M_{p,q}(\C)$, the orbit  
$K\cdot\widehat{Y}^{p,q} \subset\pgot$ is equal to the set of matrices with singular spectrum equal to $\s(Y)$. If one restricts the projection $\pi:\glgot_n(\C)\times \glgot_n(\C)\longrightarrow \glgot_n(\C)$ to the subspace $\tpgot\simeq Herm(n)$, we obtain the map 
$\pi: \tpgot\to \pgot$ that sends an Hermitian matrix $X=\begin{pmatrix}
X_{11}& X_{12}\\
X_{12}^*& X_{22}
\end{pmatrix}$ to $\widehat{X_{12}}^{p,q}=\begin{pmatrix}
0_{pp}& X_{12}\\
X_{12}^*& 0_{qq}
\end{pmatrix}$.

Since the orbit $\tK\cdot X$ is equal to $\Ocal_{\e(X)}$, we see then that $\Acal(p,q)$ can be defined as follows: $(\lambda,s)\in\R^n_+\times\R^q_{++}$ 
belongs to the cone $\Acal(p,q)$ if and only if for any $(X,Y)\in Herm(n)\times M_{p,q}(\C)$ satisfying $\lambda=\e(X)$ and $s=\s(Y)$, we have $K\cdot\widehat{Y}^{p,q}\subset \pi\left(\tK\cdot X\right)$.

%%%%%%%%%%%%%%%%%%%%%%%%%%%%%
\subsection{Inequalities determining $\Acal(p,q)$}
%%%%%%%%%%%%%%%%%%%%%%%%%%%%%

The computations done in the previous section, together with Proposition \ref{prop:oshea-sjamaar-2}, gives us the following result.

\begin{proposition} 
Let $(\lambda,s)\in\R^n_+\times\R^q_{++}$, and let $(X,Y)\in Herm(n)\times M_{p,q}(\C)$ such that $\lambda=\e(X)$ and $s=\s(Y)$.
 The following statements are equivalent:
\begin{itemize}
\item $(\lambda,s)\in \Acal(p,q)$,
\item $K\cdot\widehat{Y}^{p,q}\subset \pi\left(\tK\cdot X\right)$,
\item $U\cdot\widehat{Y}^{p,q}\subset\pi\left(\tU\cdot (X,\overline{X})\right)$,
\item $(\lambda,\lambda^*, 2\widehat{s}\,^{p,q})\in\horn(n)$.
\end{itemize}
\end{proposition}

Thanks to the description of the $\horn(n)$ cone  given in Theorem \ref{theo:horn}, we can conclude with the following description of $\Acal(p,q)$. Note that 
$|\lambda^*|_{J}=-|\lambda|_{J^o}$ and $|\widehat{s}\,^{p,q}|_{K}=|s|_{K\cap[q]}-|s|_{K^o\cap[q]}$.

\begin{proposition}\label{prop:pep=A-p-q} An element $(\lambda,s)\in\R^n_+\times \R^q_{++}$ belongs to $\Acal(p,q)$ if and only if 
$$
(\star)_{I,J,K}\qquad\qquad\qquad|\lambda|_I-|\lambda|_{J^o}\geq 2|s|_{K\cap[q]}-2|s|_{K^o\cap[q]}
$$
for any $r\leq \frac{n}{2}$ and any $(I,J,K)\in \LR^n_r$.
\end{proposition}
However, our description is less precise than that obtained by Fomin-Fulton-Li-Poon \cite{FFLP}. They show the remarkable fact that it suffices to consider inequalities 
$(\star)_{I,J,K}$ when $I,J,K$ are subsets of $[q]$.

\begin{theorem}[\cite{FFLP}]\label{theo:FFLP} An element $(\lambda,s)\in\R^n_+\times \R^q_{++}$ belongs to $\Acal(p,q)$ if and only if 
$$
|\lambda|_I-|\lambda|_{J^o}\geq 2|s|_{K}
$$
for any $r\leq q$ and any $(I,J,K)\in \LR^q_r$.
\end{theorem}

%%%%%%%%%%%%%%%%%%%%%%%%%%%%%%%%%%%%%%
\subsection{Examples}
%%%%%%%%%%%%%%%%%%%%%%%%%%%%%%%%%%%%%%

%%%%%%%%%%%%%%%%%%%%%%%%%%%%%%%%%%%%%%
\subsubsection*{Computation of $\Acal(2,2)$}\label{sec:example-2-2}
%%%%%%%%%%%%%%%%%%%%%%%%%%%%%%%%%%%%%%

The inequalities associated to $(I,J,K)\in \LR^2_1$ are 
\begin{equation}\label{eq:A-2-2-premier}
\lambda_1 - \lambda_4 \geq 2s_1,\quad \lambda_2 - \lambda_4 \geq 2s_2,\quad \lambda_1 - \lambda_3 \geq 2s_2.
\end{equation}
The inequality associated to $I=J=K=\{1,2\}$ is 
\begin{equation}\label{eq:A-2-2-second}
\lambda_1 + \lambda_2-\lambda_3-\lambda_4\geq 2(s_1+ s_2).
\end{equation}

Theorem \ref{theo:FFLP}  give us the following description.

\begin{corollary}
An element $(\lambda,s)\in \R^4_+\times\R^2_{++}$ belongs to $\Acal(2,2)$ if and only if the conditions (\ref{eq:A-2-2-premier}) and (\ref{eq:A-2-2-second}) hold.
\end{corollary}

%%%%%%%%%%%%%%%%%%%%%%%%%%%%%%%%%%%%%%
\subsubsection*{Computation of $\Acal(3,3)$}\label{sec:example-3-3}
%%%%%%%%%%%%%%%%%%%%%%%%%%%%%%%%%%%%%%

The inequalities associated to $\LR^3_1$ are 
\begin{eqnarray}\label{eq:6-1}
\lambda_1 - \lambda_6 \geq 2s_1&\quad&\lambda_1 - \lambda_4 \geq 2s_3\nonumber\\ 
 \lambda_2 - \lambda_6 \geq 2s_2& \quad& \lambda_2 - \lambda_5 \geq 2s_3 \\
 \quad \lambda_1 - \lambda_5 \geq 2s_2 & \quad& \lambda_3 - \lambda_6 \geq 2s_3.\nonumber
\end{eqnarray}
The inequalities associated to $\LR^3_2$ are 
\begin{eqnarray}\label{eq:6-2}
\lambda_1 +\lambda_2 -\lambda_5 - \lambda_6& \geq &2(s_1+s_2) \nonumber\\
\lambda_1 +\lambda_2 -\lambda_4 - \lambda_6& \geq &2(s_1+s_3) \nonumber\\
\lambda_1 +\lambda_3 -\lambda_5 - \lambda_6& \geq &2(s_1+s_3) \nonumber\\
\lambda_1 +\lambda_2 -\lambda_4 - \lambda_5& \geq &2(s_2+s_3)\\
\lambda_1 +\lambda_3 -\lambda_4 - \lambda_6& \geq &2(s_2+s_3)\nonumber \\
\lambda_2 +\lambda_3 -\lambda_5 - \lambda_6& \geq &2(s_2+s_3).\nonumber
\end{eqnarray}

The inequality associated to $I=J=K=\{1,2, 3\}$ is 
\begin{equation}\label{eq:6-3}
\lambda_1 +\lambda_2+\lambda_3 -\lambda_4-\lambda_5 - \lambda_6 \geq 2(s_1+s_2+s_3).
\end{equation}

The result of Fulton-Fomin-Li-Poon (Theorem \ref{theo:FFLP}) gives the following description of $\Acal(3,3)$.
\begin{proposition}
An element $(\lambda,s)\in \R^6_+\times\R^3_{++}$ belongs to $\Acal(3,3)$ if and only if the inequalities listed in $(\ref{eq:6-1})$, $(\ref{eq:6-2})$ and 
$(\ref{eq:6-3})$ are satisfied.
\end{proposition}

\begin{remark}The cone $\Acal(3,3)\subset \R^6\times\R^3$ corresponds to the intersection of the Horn cone 
$\horn(6)\subset \R^{18}$ with the subspace $\{(\lambda,\lambda^*,2\, \widehat{s}\,^{p,q}), (\lambda,s)\in \R^6\times \R^3\}$. Strikingly,  $\Acal(3,3)$ is determined by 
$21$ inequalities, while $\horn(6)$ is described with a minimal list of $536$ inequalities.

\end{remark}

%%%%%%%%%%%%%%%%%%%%%%%%%%%%%%%%%%%%%%
%%%%%%%%%%%%%%%%%%%%%%%%%%%%%%%%%%%%%%
\section{The cone $\Scal(p,q)$}\label{sec:S(p,q)}
%%%%%%%%%%%%%%%%%%%%%%%%%%%%%%%%%%%%%%
%%%%%%%%%%%%%%%%%%%%%%%%%%%%%%%%%%%%%%

We work with the projection $\pi_0 : \glgot_{2n}(\C)\longrightarrow\glgot_n(\C)\times \glgot_n(\C)$ defined by the relations: 
\begin{equation}\label{eq:projection-pi-0}
B= \begin{pmatrix}
B_{\bf 00}& B_{\bf 01} \\
B_{\bf 10}& B_{\bf 11} \\
\end{pmatrix}\qquad \longmapsto \qquad
\pi_0(B)= (B_{\bf 00},B_{\bf 11}).
\end{equation}
Here each matrix $B_{\bf ij}$ belongs to $\glgot_{n}(\C)$.

Recall that $(\lambda,\mu,\nu)\in \LR(n,n)$ if and only if there exists a $2n$-square Hermitian matrix $B$ such that $\lambda=\e(B)$, $\mu=\e(B_{\bf 00})$ and 
$\nu=\e(B_{\bf 11})$.

%%%%%%%%%%%%%%%%%%%%%%%%%%%%%%%%%%%%%%%%%%%%%%%
\subsection{Matrix identities}\label{sec:matrix-S-p-q}
%%%%%%%%%%%%%%%%%%%%%%%%%%%%%%%%%%%%%%%%%%%%%%%%%%
Here we use the notations $\widehat{Y}^{p,q}, \widehat{\mu}^{p,q}$  introduced in the \S \ref{sec:augmented}.

Let us decompose a $n$-square complex matrix 
$X=\begin{pmatrix}
X_{11}& X_{12}\\
X_{21}& X_{22}
\end{pmatrix}$
by blocks where $X_{12}\in M_{p,q}(\C)$.  Let 
$\widehat{X}= \begin{pmatrix}
0& X \\
X^*& 0 \\
\end{pmatrix}$
be the associated $2n$-square Hermitian matrix. Recall that  $\e(\widehat{X})=\widehat{\s(X)}$ (see Section \ref{sec:augmented}).

Let $P_\tau\in O_{2n}(\R)$ be the orthogonal matrix associated with the permutation $\tau:[2n]\to[2n]$ which is defined as follows: $\tau(k)=k$ if $1\leq k\leq p$, $\tau(k)=k+q$ if $p+1\leq k\leq n+p$ and 
$\tau(k)=k-n$ if $n+p+1\leq k\leq 2n$.

We see then that $P_\tau \widehat{X} P_\tau^{-1}$ is a $2n$-square hermitian matrix such that 
$$
(P_\tau \widehat{X} P_\tau^{-1})_{\bf 00}=
\begin{pmatrix}
0& X_{12} \\
X_{12}^*& 0 \\
\end{pmatrix}
\quad {\rm and}\quad 
(P_\tau \widehat{X} P_\tau^{-1})_{\bf 11}=
\begin{pmatrix}
0& X_{21} \\
X_{21}^*& 0 \\
\end{pmatrix}
$$
Finally we obtain the relations $\widehat{\s(X)}=\e(\widehat{X})=\e(P_\tau \widehat{X} P_\tau^{-1})$,  
$$
\e((P_\tau \widehat{X} P_\tau^{-1})_{\bf 00})=\widehat{\s(X_{12})}^{p,q}\quad {\rm and}\quad \e((P_\tau \widehat{X} P_\tau^{-1})_{\bf 11})=\widehat{\s(X_{21})}^{p,q}.
$$
In other words, $(\widehat{\s(X)},\widehat{\s(X_{12})}^{p,q},\widehat{\s(X_{21})}^{p,q})\in \LR(n,n)$ for any  $n$-square complex matrix $X$. At this point, we have shown 
that any $(\gamma, s,t)\in \Scal(p,q)$ satisfies the relation
\begin{equation}\label{eq:OSS-relation-S-p-q}
(\widehat{\gamma},\widehat{s}\,^{p,q}, \widehat{t}\,^{p,q})\in \LR(n,n).
\end{equation}

In the next sections, we explain how the O'Shea-Sjamaar theorem (see Proposition \ref{prop:oshea-sjamaar-2}) allows us to show that (\ref{eq:OSS-relation-S-p-q}) characterizes the cone $\Scal(p,q)$.

%%%%%%%%%%%%%%%%%%%%%%%%%%%%%%%%%%%%%%%%%%%%%%%
\subsection{Antiholomorphic involution and orthogonal projection}
%%%%%%%%%%%%%%%%%%%%%%%%%%%%%%%%%%%%%%%%%%%%%%%%%%

We work with the real reductive Lie groups $G:= U(p,q)\times U(q,p)$ and $\tG:= U(n,n)$. The embedding $\iota: G\to \tG$ is defined as follows:
\begin{equation}\label{eq:embedding-p-q}
\iota(g,h)= \begin{pmatrix}
g_{11}& 0_{pn} & g_{12}\\
0_{np}& h & 0_{nq}\\
g_{21}& 0_{qn} & g_{22}
\end{pmatrix},\qquad{\rm when}\qquad g= \begin{pmatrix}
g_{11}&  g_{12}\\
g_{21}&  g_{22}
\end{pmatrix}.
\end{equation}
Here $g_{11}\in M_{p,p}(\C)$, $g_{12}\in M_{p,q}(\C)$, $g_{2,1}\in M_{q,p}(\C)$ and $g_{22}\in M_{q,q}(\C)$.

The unitary group $\tK:=U_n\times U_n$ is a maximal compact subgroup of $\tG$. The subspace 
$\tpgot:= \{\widehat{X}, X\in \glgot_{n}(\C)\}\subset\tggot$ 
admits a canonical action of $\tK$. The subgroup $K= K_1\times K_2$, with $K_1\simeq U_p\times U_q$ and $K_2\simeq U_q\times U_p$, is a maximal compact subgroup of $G$, 
and the subspace $\pgot=\tpgot\cap\ggot$ admits a natural identification with $M_{p,q}(\C)\times M_{q,p}(\C)$:
$$
(Y,Z)\in M_{p,q}(\C)\times M_{q,p}(\C)\longmapsto (\widehat{Y}^{p,q},\widehat{Z}^{q,p})\in \pgot.
$$

The complexification of the group $G$ is $G_\C:= GL_n(\C)\times GL_n(\C) $. We consider the antiholomorphic involution $\sigma$ on $G_\C$ defined by 
$\sigma(g,h)= (I_{p,q}(g^*)^{-1}I_{p,q},I_{q,p}(h^*)^{-1}I_{q,p})$. The subgroup $G$ is the fixed point set of $\sigma$.

The complexification of the group $\tG$ is $\tG_\C:= GL_{2n}(\C)$. We consider the antiholomorphic involution $\widetilde{\sigma}$ on $\tG_\C$ defined by 
$\widetilde{\sigma}(g)= I_{n,n}(g^*)^{-1}I_{n,n}$. 
The subgroup $\tG$ corresponds to the fixed point set of $\widetilde{\sigma}$.

The groups $U=U_n\times U_n$ and $\tU= U_{2n}$  are respectively  maximal compact subgroups of $G_\C$ and $\tG_\C$. The fixed point subgroups of the involutions are $U^\sigma=K$ and $\tU^{\tilde{\sigma}}=\tK$.

The embedding $\iota: G\croc \tG$ admits a complexification $\iota_\C: G_\C\croc \tG_\C$. At the level of Lie algebra, we have a morphism 
$\iota_\C: \glgot_n(\C)\times \glgot_n(\C)\croc \glgot_{2n}(\C)$, still defined by (\ref{eq:embedding-p-q}).

The orthogonal projection $\pi_1 : \glgot_{2n}(\C)\longrightarrow\glgot_n(\C)\times \glgot_n(\C)$ dual to the morphism $\iota_\C$ 
is defined by the relations : 
\begin{equation}\label{eq:projection-p-q}
A= \begin{pmatrix}
A_{11}& A_{12} & A_{13}\\
A_{21}& A_{22} & A_{23}\\
A_{31}& A_{32} & A_{33}\\
\end{pmatrix}\qquad \longmapsto \qquad\pi_1(A)= \left(\begin{pmatrix}
A_{11}&  A_{13}\\
A_{31}&  A_{33}
\end{pmatrix}, A_{22}\right).
\end{equation}
Here the matrix $A\in\glgot_{2n}(\C)$ is written by blocks relatively to the decomposition $2n=p+n+q$.

In the beginning of \S \ref{sec:S(p,q)}, we have consider another projection $\pi_0$ (see (\ref{eq:projection-pi-0})).

\begin{lemma}\label{lem:projection-pi-0}
For any $U_{2n}$-orbit $\Ocal \subset \glgot_{2n}(\C)$, we have $\pi_1(\Ocal)=\pi_0(\Ocal)$.
\end{lemma}

{\em Proof:} Let $P_\tau\in O_{2n}(\R)$ the orthogonal matrix defined in Section \ref{sec:matrix-S-p-q}. We check that $\pi_1(M)=\pi_0(P_\tau M P^{-1}_\tau)$, $\forall M\in\glgot_{2n}(\C)$. Our lemma follows from this relation. $\Box$

%%%%%%%%%%%%%%%%%%%%%%%%%%%%%%%%
\subsection{Description of $\Scal(p,q)$ through $\LR(n,n)$}
%%%%%%%%%%%%%%%%%%%%%%%%%%%%%%%%

For $t\in \R^q$, we consider the $n$-square Hermitian matrix 

\begin{equation}\label{eq:X-c}
Y(t):=\begin{pmatrix}0_{qq}&  0_{q,p-q}& \diag(t)\\
0_{p-q,q}&  0_{p-q,p-q}& 0_{p-q,q}\\
\diag(t) & 0_{q,p-q} & 0_{qq}
\end{pmatrix}.
\end{equation}

Here is the main application of Proposition \ref{prop:oshea-sjamaar-2}.

\begin{proposition}\label{prop:equivalence-A-p-q}
Let $(\gamma,s,t)\in \R^n_{++}\times\R^q_{++}\times\R^q_{++}$. The following statements are equivalent:
\begin{enumerate}
\item $(\gamma,s,t)\in\Scal(p,q)$,
\item $\exists A= \begin{pmatrix}
A_{11}& A_{12} \\
A_{21}& A_{22} \\
\end{pmatrix}\in \glgot_{n}(\C)$, such that $\s(A)=\gamma$, $\s(A_{12})=s$, and $\s(A_{21})=t$,
\item $\exists M\in\tpgot$, with $\pi_1(M)=(M_1,M_2)$, s.t. $\e(M)=\widehat{\gamma}$, $\e(M_1)=\widehat{s}\,^{p,q}$ and $\e(M_2)=\widehat{t}\,^{p,q}$,
\item $\pi_1\left(U_n\times U_n\cdot \widehat{\diag(\gamma)}\right)$ contains 
$\Big(U_p\times U_q\cdot Y(s)\Big)\times \Big(U_q\times U_p\cdot Y(t)\Big)$.
\item $\pi_1\left(U_{2n}\cdot \widehat{\diag(\gamma)}\right)$ contains 
$U_n\cdot Y(s)\times U_n\cdot Y(t)$.
\item $\pi_0\left(U_{2n}\cdot \diag(\widehat{\gamma})\right)$ contains 
$\Ocal_{\widehat{s}\,^{p,q}}\times \Ocal_{\widehat{t}\,^{p,q}}$.
\item $(\widehat{\gamma},\widehat{s}\,^{p,q}, \widehat{t}\,^{p,q})\in \LR(n,n)$.
\end{enumerate}
\end{proposition}

{\em Proof:}  Equivalences ``$1. \Longleftrightarrow 2.$''  and ``$6. \Longleftrightarrow 7.$'' are true by definition.  Equivalence ``$2. \Longleftrightarrow 3.$'' is proved by taking $M=\widehat{A}$ (see \S \ref{sec:matrix-S-p-q}).
Equivalence ``$3. \Longleftrightarrow 4.$'' is obtained by noting the following relations 
$$
\Big\{M\in\tpgot;\ \e(M)=\widehat{\gamma}\Big\}=U_n\times U_n\cdot \widehat{\diag(\gamma)},
$$
and $\big\{(X,Y)\in\pgot;\ \e(X)=\widehat{s}\ {\rm and} \ \e(Y)=\widehat{t}\,\big\}= \left(U_p\times U_q\cdot Y(s)\right)\times \left(U_q\times U_p\cdot Y(t)\right)$.
Equivalence ``$4. \Longleftrightarrow 5.$'' follows from Proposition \ref{prop:oshea-sjamaar-2}, and ``$5. \Longleftrightarrow 6.$''
is a consequence of Lemma \ref{lem:projection-pi-0} and the fact that the orbit $U_{2n}\cdot \widehat{\diag(\gamma)}$ is equal to 
$U_{2n}\cdot \diag(\widehat{\gamma})$. $\Box$

%%%%%%%%%%%%%%%%%%%%%%%%%%%%%
\subsection{Inequalities determining $\Scal(p,q)$}
%%%%%%%%%%%%%%%%%%%%%%%%%%%%%

Thanks to Proposition \ref{prop:equivalence-A-p-q} and Theorem \ref{theo:LR-m-n}, we obtain the following description of the cone $\Scal(p,q)$.

\begin{theorem}\label{theo:S-p-q} An element $(\gamma,s,t)\in \R^n_{++}\times\R^q_{++}\times\R^q_{++}$ belongs to $\Scal(p,q)$ if and only if, 
\begin{equation}\label{eq:S-p-q-majoration}
\gamma_k\geq s_k\quad  {\rm and}\quad  \gamma_k\geq t_k,\qquad \forall k\in [q],
\end{equation}
and 
\begin{equation}\label{eq:theoreme-S-p-q}
|\gamma |_{A\cap[n]} \, -\,  |\gamma |_{A^o\cap[n]}\ \geq\   | s |_{B\cap[q]} \,-\,  | s |_{B^o\cap[q]} \,+\,   | t |_{C\cap[q]} \,-\,  | t |_{C^o\cap[q]},
\end{equation}
holds for any triplets $(A,B,C)$  satisfying the following conditions : 
\begin{itemize}
\item $B,C$ are strict subsets of $[n]$,
\item $A\subset [2n]$ and $\sharp A = \sharp B + \sharp C$,
\item the Littlewood-Richardson coefficient $c^{\mu(A)}_{\mu(B),\mu(C)}$ is non-zero.
\end{itemize}
\end{theorem}

Let us use some duality to minimize the number of equations (see (\ref{eq:duality-equations})). The equation (\ref{eq:theoreme-S-p-q}) means that  
$(\widehat{\gamma},\widehat{s}\,^{p,q},\widehat{t}\,^{p,q})$ satisfies 
$|x|_A\geq |y|_B + |z|_C$, that is $|\widehat{\gamma}|_A\geq |\widehat{s}\,^{p,q} |_B + |\widehat{t}\,^{p,q}|_C$. 
If we apply $(\widehat{\gamma},\widehat{s}\,^{p,q},\widehat{t}\,^{p,q})$ to the relation $|x|_{A^o}\leq |y|_{C^o} + |z|_{C^o}$, we get 
$|\widehat{\gamma}|_{A^o}\leq |\widehat{s}\,^{p,q} |_{B^o} + |\widehat{t}\,^{p,q}|_{C^o}$ which is equivalent to (\ref{eq:theoreme-S-p-q}) since  
$|\widehat{\gamma}|_{A^o}=-|\widehat{\gamma}|_{A}$,  $|\widehat{s}\,^{p,q} |_{B^o}=-|\widehat{s}\,^{p,q} |_{B}$ and  
$|\widehat{t}\,^{p,q}|_{C^o}=-|\widehat{t}\,^{p,q}|_{C}$.

\medskip

We can therefore rewrite Theorem \ref{theo:S-p-q}, by requiring that (\ref{eq:theoreme-S-p-q}) holds for all strict subsets $A\subset [2n]$, $B,C\subset [n]$, which satisfy $\sharp A =\sharp B +\sharp C\leq n$ and $c^{\mu(A)}_{\mu(B),\mu(C)}\neq 0$.

%%%%%%%%%%%%%%%%%%%%%%%%%%%%%%%%%%%%
\subsection{Examples}
%%%%%%%%%%%%%%%%%%%%%%%%%%%%%%%%%%%%

\subsubsection*{The cone $\Scal(1,1)$}

We have to look to subsets $B=\{b\}\subset [2]$, $C=\{c\}\subset [2]$, and $A=\{a_2>a_1\}\subset [4]$ such that the Littlewood-Richardson coefficient $c^{\mu(A)}_{\mu(B),\mu(C)}$ is non-zero. Here's the list and the corresponding inequalities:
\begin{itemize}
\item[$i)$] $B=C=\{2\}$ and $A=\{4,1\}$ or $\{3,2\}$: \quad $0\geq -s-t$.
\item[$ii)$] $B=\{1\}$, $C=\{2\}$ and $A=\{3,1\}$: \quad $\lambda_1-\lambda_2\geq s-t$.
\item[$iii)$] $B=\{2\}$, $C=\{1\}$ and $A=\{3,1\}$: \quad $\lambda_1-\lambda_2\geq -s+t$.
\item[$iv)$] $B=C=\{1\}$ and $A=\{2,1\}$: \quad $\lambda_1+\lambda_2\geq s+t$.
\end{itemize}

Note that the inequalities (\ref{eq:S-p-q-majoration}) are here a consequence of $ii)$, $iii)$ and $iv)$. Thus, an element $(\gamma,s,t)\in \R^2_{++}\times\R^{\geq 0}\times\R^{\geq 0}$ 
belongs to $\Scal(1,1)$ if and only if
$$
\lambda_1-\lambda_2\geq |s-t|\quad {\rm and}\quad \lambda_1+\lambda_2\geq s+t.
$$
We recover the computation done in \cite{FFLP} (see Example 1.17).

\subsubsection*{The cone $\Scal(2,1)$}

First, we look to subsets $B=\{b\}\subset [3]$, $C=\{c\}\subset [3]$, and $A=\{a_2>a_1\}\subset [6]$ such that the Littlewood-Richardson coefficient $c^{\mu(A)}_{\mu(B),\mu(C)}$ is non-zero. 
The corresponding inequality (\ref{eq:theoreme-S-p-q}) is called trivial when it is a consequence 
of the following relations
\begin{equation}\label{eq:inequalities-S-2-1-first}
\gamma_1\geq\gamma_2\geq \gamma_3\geq 0,\qquad \gamma_1\geq t\geq 0,\qquad  \gamma_1\geq s\geq 0.
\end{equation}
Here's the list of the non-trivial inequalities:
\begin{itemize}
\item $B=\{1\}$, $C=\{1\}$ and $A=\{2,1\}$: \quad $ \gamma_1+ \gamma_2\geq s+t$.
\item $B=\{1\}$, $C=\{3\}$ and $A=\{4,1\}$: \quad $ \gamma_1 -  \gamma_3\geq s-t$.
\item $B=\{3\}$, $C=\{1\}$ and $A=\{4,1\}$: \quad $ \gamma_1 -  \gamma_3\geq -s+t$.
%\item $B=\{1\}$, $C=\{2\}$ and $A=\{3,1\}$: \quad $\lambda_1+\lambda_3\geq s$.
%\item $B=\{2\}$, $C=\{1\}$ and $A=\{3,1\}$: \quad $\lambda_1+\lambda_3\geq t$.
\end{itemize}

Next, we examine subsets $B=\{b_1>b_2\}\subset [3]$, $C=\{c\}\subset [3]$, and $A=\{a_3>a_2>a_1\}\subset [6]$ such that $c^{\mu(A)}_{\mu(B),\mu(C)}\neq 0$. An easy check 
shows that all inequalities obtained here are a consequence of (\ref{eq:inequalities-S-2-1-first}) and the inequalities
\begin{equation}\label{eq:inequalities-S-2-1-second}
\lambda_1+\lambda_2\geq s+t,\qquad \lambda_1-\lambda_3\geq |s-t|,
\end{equation}
which we have just proved above.

\begin{corollary}
An element $(\gamma,s,t)\in \R^3\times\R\times\R$ belongs to $\Scal(2,1)$ if and only if the inequalities (\ref{eq:inequalities-S-2-1-first}) and (\ref{eq:inequalities-S-2-1-second}) hold.
\end{corollary}

%%%%%%%%%%%%%%%%%%%%%%%%%%%%%%%%%%%%%%
%%%%%%%%%%%%%%%%%%%%%%%%%%%%%%%%%%%%%%
\section{The cone $\Tcal(p,q)$}\label{sec:T(p,q)}
%%%%%%%%%%%%%%%%%%%%%%%%%%%%%%%%%%%%%%
%%%%%%%%%%%%%%%%%%%%%%%%%%%%%%%%%%%%%%

We consider here the projections $\pi_0,\pi_1 : \glgot_{2n}(\C)\longrightarrow\glgot_{2p}(\C)\times \glgot_{2q}(\C)$:  
\begin{itemize}
\item $\pi_1$ is defined by the relations: 
\begin{equation}\label{eq:projection-T-p-q}
A= \begin{pmatrix}
A_{11}& A_{12} & A_{13}\\
A_{21}& A_{22} & A_{23}\\
A_{31}& A_{32} & A_{33}\\
\end{pmatrix}\qquad \longmapsto \qquad\pi(A)= \left(\begin{pmatrix}
A_{11}&  A_{13}\\
A_{31}&  A_{33}
\end{pmatrix}, A_{22}\right),
\end{equation}
where the matrix $A\in\glgot_{2n}(\C)$ is written by blocks relatively to the decomposition $2n=p+2q+p$.
\item $\pi_0$ is defined by the relations: 
\begin{equation}\label{eq:projection-T-p-q-bis}
B= \begin{pmatrix}
B_{\bf 00}&B_{\bf 01}\\
B_{\bf 10}&B_{\bf 11}\\
\end{pmatrix}\qquad \longmapsto \qquad
\pi_0(B)= (B_{\bf 00},B_{\bf 11}),
\end{equation}
where $B_{\bf 00}\in \glgot_{2p}(\C)$ and $B_{\bf 11}\in \glgot_{2q}(\C)$
\end{itemize}

\begin{lemma}\label{lem:projection-pi-0-bis}
For any $U_{2n}$-orbit $\Ocal \subset \glgot_{2n}(\C)$, we have $\pi_1(\Ocal)=\pi_0(\Ocal)$.
\end{lemma}

{\em Proof:} Same proof as for Lemma \ref{lem:projection-pi-0}. $\Box$

%%%%%%%%%%%%%%%%%%%%%%%%%%%%%%%%%%%%%%%%%%%%%%%
\subsection{Matrix identities}\label{sec:matrix-T-p-q}
%%%%%%%%%%%%%%%%%%%%%%%%%%%%%%%%%%%%%%%%%%%%%%%

Let us decompose a $n$-square complex matrix 
$X=\begin{pmatrix}
X_{11}& X_{12}\\
X_{21}& X_{22}
\end{pmatrix}$
by blocks where $X_{12}\in M_{p,q}(\C)$.  We want to find a link between the singular eigenvalues of $X$, $X_{11}$ and $X_{22}$.

The matrix $Q=\begin{pmatrix}
0_{pq}& Id_p \\
Id_q& 0_{qp} \\
\end{pmatrix}$ is orthogonal  and the matrix 
$X':=XQ= \begin{pmatrix}
X_{12}& X_{11}\\
X_{22}& X_{21}
\end{pmatrix}$ has the same singular values as $X$. The image of the $2n$-square Hermitian matrix 
$\widehat{X'}:=\begin{pmatrix}
0_{nn}& X'\\
(X')^*& 0_{nn}
\end{pmatrix}$ trough the projection $\pi_1$ is equal to 
$$
\pi_1\left(\widehat{X'}\right)=\left(
\widehat{X_{11}},\widehat{X_{22}}
\right).
$$
The orbit $\Ocal:=U_{2n}\cdot\widehat{X'}$ is equal to the subset of $2n$-square Hermitian matrices $Y$ satisfying $\e(Y)=\e\left(\widehat{X'}\right)=\widehat{\s(X')}=\widehat{\s(X)}$, and 
the projection $\pi_0(\Ocal)=\pi_1(\Ocal)$ contains $(\widehat{X_{11}},\widehat{X_{22}})$, so $(\widehat{\s(X)}, \widehat{\s(X_{11})}, \widehat{\s(X_{22})})\in \LR(2p,2q)$.

We have just shown that any $(\gamma, s,t)\in \Tcal(p,q)$ satisfies the relation
\begin{equation}\label{eq:OSS-relation-T-p-q}
(\,\widehat{\gamma},\widehat{s}, \widehat{t}\,)\in \LR(2p,2q).
\end{equation}

In the next sections, we explain how the O'Shea-Sjamaar theorem (see Proposition \ref{prop:oshea-sjamaar-2}) allows us to see that relation (\ref{eq:OSS-relation-T-p-q}) 
characterizes the cone $\Tcal(p,q)$.

%%%%%%%%%%%%%%%%%%%%%%%%%%%%%%%%%%%%%%%%%%%%%%%
\subsection{Antiholomorphic involution and orthogonal projection}
%%%%%%%%%%%%%%%%%%%%%%%%%%%%%%%%%%%%%%%%%%%%%%%%%%

We work with the real reductive Lie groups $G:= U(p,p)\times U(q,q)$ and $\tG:= U(n,n)$. The embedding $\iota: G\to \tG$ is defined as follows:
\begin{equation}\label{eq:embedding-T-p-q}
\iota(g,h)= \begin{pmatrix}
g_{11}& 0_{p,2q} & g_{12}\\
0_{2q,p}& h & 0_{2q,p}\\
g_{21}& 0_{p,2q} & g_{22}
\end{pmatrix},\qquad{\rm when}\qquad g= \begin{pmatrix}
g_{11}&  g_{12}\\
g_{21}&  g_{22}
\end{pmatrix}.
\end{equation}
Here $g_{ij}\in \glgot_{p}(\C)$ and $h\in U(q,q)\subset \glgot_{2q}(\C)$.

The unitary group $\tK:=U_n\times U_n$ is a maximal compact subgroup of $\tG$. The subspace 
$\tpgot:= \{\widehat{X}, X\in M_{n,n}(\C)\}\subset\tggot$ 
admits a canonical action of $\tK$. The subgroup $K= K_1\times K_2$, with $K_1\simeq U_p\times U_p$ and $K_2\simeq U_q\times U_q$, is a maximal compact subgroup of $G$, 
and the subspace $\pgot=\tpgot\cap\ggot$ admits a natural identification with $\glgot_{p}(\C)\times \glgot_{q}(\C)$:
$$
(Y,Z)\in \glgot_{p}(\C)\times \glgot_{q}(\C)\longmapsto (\widehat{Y},\widehat{Z})\in \pgot.
$$
%See \S \ref{sec:augmented} for the notations $\widehat{X}$.
%$\pgot=\pgot_1\times\pgot_2$

The complexification of the group $G$ is $G_\C:= GL_{2p}(\C)\times GL_{2q}(\C) $. We consider the antiholomorphic involution $\sigma$ on $G_\C$ defined by 
$\sigma(g,h)= (I_{p,p}(g^*)^{-1}I_{p,p},I_{q,q}(h^*)^{-1}I_{q,q})$. The subgroup $G$ is the fixed point set of $\sigma$.

The complexification of the group $\tG$ is $\tG_\C:= GL_{2n}(\C)$. We consider the antiholomorphic involution $\widetilde{\sigma}$ on $\tG_\C$ defined by 
$\widetilde{\sigma}(g)= I_{n,n}(g^*)^{-1}I_{n,n}$. 
The subgroup $\tG$ corresponds to the fixed point set of $\widetilde{\sigma}$.

The groups $U=U_{2p}\times U_{2q}$ and $\tU= U_{2n}$  are respectively  maximal compact subgroups of $G_\C$ and $\tG_\C$. The fixed point subgroups of the involutions are $U^\sigma=K$ and $\tU^{\tilde{\sigma}}=\tK$.

The embedding $\iota: G\croc \tG$ admits a complexification $\iota_\C: G_\C\croc \tG_\C$. At the level of Lie algebra, we have a morphism 
$\iota_\C: \glgot_{2p}(\C)\times \glgot_{2q}(\C)\croc \glgot_{2n}(\C)$, still defined by (\ref{eq:embedding-T-p-q}).

The orthogonal projection $\glgot_{2n}(\C)\longrightarrow\glgot_{2p}(\C)\times \glgot_{2q}(\C)$ dual to the morphism $\iota_\C$ is the map $\pi_1$ defined at the start of 
Section \ref{sec:T(p,q)}.

%%%%%%%%%%%%%%%%%%%%%%%%%%%%%%%%
\subsection{Description of $\Tcal(p,q)$ through $\LR(2p,2q)$}
%%%%%%%%%%%%%%%%%%%%%%%%%%%%%%%%

Here is the main application of the Proposition \ref{prop:oshea-sjamaar-2}. Recall that for $(\gamma,s,t)\in \R^n_{++}\times\R^p_{++}\times\R^q_{++}$, we define
$\widehat{\gamma}:=(\gamma_1,\ldots,\gamma_n,-\gamma_n,\ldots,-\gamma_1)$, $\widehat{s}:=(s_1,\ldots,s_p,-s_p,\ldots,-s_1)$ and 
$\widehat{t}:=(t_1,\ldots,t_q,-t_q,\ldots,-t_1)$.

\begin{proposition}\label{prop:equivalence-T-p-q}
Let $(\gamma,s,t)\in \R^n_{++}\times\R^p_{++}\times\R^q_{++}$. The following statements are equivalent:
\begin{enumerate}
\item $(\gamma,s,t)\in\Tcal(p,q)$,
\item $\exists A= \begin{pmatrix}
A_{11}& A_{12} \\
A_{21}& A_{22} \\
\end{pmatrix}\in \glgot_{n}(\C)$, such that $\s(A)=\gamma$, $\s(A_{11})=s$, and $\s(A_{22})=t$,
\item $\exists M\in\tpgot$, with $\pi_1(M)=(M_1,M_2)$, s.t. $\e(M)=\widehat{\gamma}$, $\e(M_1)=\widehat{s}$ and $\e(M_2)=\widehat{t}$,
\item $\pi_1\left(U_n\times U_n\cdot \widehat{\diag(\gamma)}\right)$ contains 
$\left(U_p\times U_p\cdot \widehat{\diag(s)}\right)\times \left(U_q\times U_q\cdot \widehat{\diag(t)}\right)$.
\item $\pi_1\left(U_{2n}\cdot \widehat{\diag(\gamma)}\right)$ contains 
$U_{2p}\cdot \widehat{\diag(s)}\times U_{2q}\cdot \widehat{\diag(t)}$.
\item $\pi_0\left(U_{2n}\cdot \diag(\widehat{\gamma})\right)$ contains 
$U_{2p}\cdot \diag(\widehat{s}\,)\times U_{2q}\cdot \diag(\widehat{t}\,)$.
\item $(\,\widehat{\gamma}\,,\,\widehat{s}\,,\, \widehat{t}\,)\in \LR(2p,2q)$.
\end{enumerate}
\end{proposition}

{\em Proof:} Equivalences ``$2. \Longleftrightarrow 3.$'' and ``$6. \Longleftrightarrow 7.$'' are true by definition.  Equivalence ``$2. \Longleftrightarrow 3.$'' is proved by taking $M=\widehat{A'}$ (see \S \ref{sec:matrix-T-p-q}). Equivalence ``$3. \Longleftrightarrow 4.$'' is obtained by noting the following relations 
$$
\Big\{M\in\tpgot;\ \e(M)=\widehat{\gamma}\Big\}=U_n\times U_n\cdot \widehat{\diag(\gamma)},
$$
and $\big\{(X,Y)\in\pgot;\ \e(X)=\widehat{s}\ {\rm and} \ \e(Y)=\widehat{t}\,\big\}= \left(U_p\times U_p\cdot \widehat{\diag(s)}\right)\times \left(U_q\times U_q\cdot \widehat{\diag(t)}\right)$.
Equivalence ``$4. \Longleftrightarrow 5.$'' follows from Proposition \ref{prop:oshea-sjamaar-2}, and ``$5. \Longleftrightarrow 6.$'' 
is a consequence of Lemma \ref{lem:projection-pi-0}. $\Box$

%%%%%%%%%%%%%%%%%%%%%%%%%%%%%
\subsection{Inequalities determining $\Tcal(p,q)$}
%%%%%%%%%%%%%%%%%%%%%%%%%%%%%

Thanks to Proposition \ref{prop:equivalence-T-p-q} and Theorem \ref{theo:LR-m-n}, we obtain the following description of the cone $\Tcal(p,q)$.

\begin{theorem}\label{theo:T-p-q} An element $(\gamma,s,t)\in \R^n_{++}\times\R^p_{++}\times\R^q_{++}$ belongs to $\Tcal(p,q)$ if and only if the following inequalities hold:
\begin{enumerate}
\item $\gamma_k\geq s_k$, $\forall k\leq p$,
\item $\gamma_j\geq t_j$, $\forall j\leq q$, 
\item $\gamma_{2q+\ell}\leq s_\ell$, $\forall \ell\leq p-q$,
\item we have $|\gamma |_{A\cap[n]} \, -\,  |\gamma |_{A^o\cap[n]}\ \geq\   | s |_{B\cap[p]} \,-\,  | s |_{B^o\cap[p]} \,+\,   | t |_{C\cap[q]} \,-\,  | t |_{C^o\cap[q]}$
\medskip

for any triplets $(A,B,C)$  satisfying the following conditions : 
\begin{itemize}
\item $B\subset [2p]$ and $C\subset [2q]$ are strict subsets,
\item $A\subset [2n]$ and $\sharp A = \sharp B + \sharp C\leq n$,
\item the Littlewood-Richardson coefficient $c^{\mu(A)}_{\mu(B),\mu(C)}$ is non-zero.
\end{itemize}
\end{enumerate}
\end{theorem}

%%%%%%%%%%%%%%%%%%%%%%%%%%%%%%%%%%%%
\subsection{Examples}
%%%%%%%%%%%%%%%%%%%%%%%%%%%%%%%%%%%%

\subsubsection*{The cone $\Tcal(1,1)$}

We need to find subsets $B=\{b\}\subset [2]$, $C=\{c\}\subset [2]$, and $A=\{a_2>a_1\}\subset [4]$ such that the Littlewood-Richardson coefficient $c^{\mu(A)}_{\mu(B),\mu(C)}$ is non-zero. This work has been done for the cone $\Scal(1,1)$. Thus, an element $(\gamma,s,t)\in \R^2_{++}\times\R^{\geq 0}\times\R^{\geq 0}$ 
belongs to $\Tcal(1,1)$ if and only if
\begin{equation}\label{eq:T-1-1}
\gamma_1-\gamma_2\geq |s-t|\quad {\rm and}\quad \gamma_1+\gamma_2\geq s+t.
\end{equation}
Note that inequalities $\gamma_1\geq s$ and $\gamma_1\geq t$ follow from (\ref{eq:T-1-1}).

\subsubsection*{The cone $\Tcal(2,1)$}

We work with $(\gamma,s,t)\in \R^3_{++}\times\R^2_{++}\times\R^{\geq 0}$ satisfying
\begin{equation}\label{eq:T-2-1-A}
\gamma_1\geq s_1, \quad \gamma_1\geq t,\quad \gamma_2\geq s_2,\quad \gamma_{3}\leq s_1.
\end{equation}

We're now interested in the inequalities associated with triplets $(A,B,C)$ such that $c^{\mu(A)}_{\mu(B),\mu(C)}$ is non-zero and $\sharp A = \sharp B + \sharp C\leq 3$. We obtain the following inequalities when $\sharp A=2$:
\begin{eqnarray}\label{eq:T-2-1-B}
\gamma_1+\gamma_2&\geq& s_1+t\nonumber\\ 
\gamma_1+\gamma_3&\geq& s_2+t \nonumber\\ 
\gamma_1-\gamma_3&\geq& |s_2-t| \\ 
\gamma_1-\gamma_2&\geq& -s_1+t.\nonumber
\end{eqnarray}
In the list (\ref{eq:T-2-1-B}), I haven't included inequality $\gamma_2+\gamma_3\geq s_2-t$, which is associated with $A=\{2,3\}$, $B=\{2\}$, $C=\{2\}$, since it follows from (\ref{eq:T-2-1-A}) and the fact that $\gamma_3,t\geq 0$.
When $\sharp A=3$ we get:
\begin{eqnarray}\label{eq:T-2-1-C}
\gamma_1+\gamma_2+\gamma_3&\geq& s_1+s_2+t\nonumber\\ 
\gamma_1+\gamma_2-\gamma_3&\geq& s_1-s_2+t \nonumber\\ 
\gamma_1+\gamma_2-\gamma_3&\geq& s_1+s_2-t \\ 
\gamma_1-\gamma_2+\gamma_3&\geq& -s_1+s_2+t \nonumber\\ 
\gamma_1-\gamma_2-\gamma_3&\geq& -s_1+s_2-t \nonumber\\ 
\gamma_1-\gamma_2-\gamma_3&\geq& -s_1-s_2+t.  \nonumber
\end{eqnarray}

Thus, an element $(\gamma,s,t)\in \R^3_{++}\times\R^2_{++}\times\R^{\geq 0}$
belongs to $\Tcal(2,1)$ if and only if the inequalities (\ref{eq:T-2-1-A}), (\ref{eq:T-2-1-B}) and (\ref{eq:T-2-1-C}) are satisfied.

%%%%%%%%%%%%%%%%%%%%%%%%%%%%%
\subsection{Interlacing inequalities for singular values}
%%%%%%%%%%%%%%%%%%%%%%%%%%%%%

Let us consider the case where $p\geq q=1$. 

Let $\gamma_1\geq\cdots\geq\gamma_{p+1}\geq 0$ be the singular values of a $p+1$-square complex matrix $X$. 
Let $X'$ be the $p$-square submatrix of $X$ obtained by deleting a row and a column: we denote by 
$s_1\geq\cdots\geq s_{p}\geq 0$ its singular spectrum.

Points 1. and 3. of Theorem \ref{theo:T-p-q} yields interlacing inequalities which where first observed by Thompson \cite{Thompson72}:
\begin{align*} 
\gamma_3  \leq s_1 & \leq  \gamma_1, \\ 
\gamma_4 \leq s_2 & \leq  \gamma_2, \\ 
& \cdots \\ 
\gamma_{j+2} \leq s_j & \leq  \gamma_j, \qquad 1\leq j\leq p-2, \\ 
& \cdots \\ 
\gamma_p  \leq s_{p-2} & \leq  \gamma_{p-2}, \\ 
s_{p-1} & \leq  \gamma_{p-1}, \\ 
s_{p} & \leq  \gamma_{p}. \\ 
\end{align*}

\end{document}